\documentclass[reqno]{amsart}
\usepackage{amsmath}
\usepackage{amssymb}
\usepackage{xspace}
\usepackage{epic}

\makeatletter

\renewcommand{\subsection}{%
\setcounter{subsection}{\value{theorem}}%
\stepcounter{theorem}%
\@startsection
{subsection}%
{1}% level
{0em}% indent
{-\baselineskip}% reduce space above
{-4pt}% negative makes a run-in heading, space before heading text
{\bfseries\normalsize}}%style
\renewcommand{\p@enumi}{\thetheorem.}
\renewcommand{\p@enumii}{\thetheorem.\theenumi.}
\makeatother
\newcommand{\Rev}[1]{\textit{Rev}(#1)}
\newcommand{\Weld}[1]{\textit{Weld}(#1)}
\newcommand{\Rules}{\textit{Rules}\xspace}
\newcommand{\SLtwo}{\textit{SL2}}
\newcommand{\NRed}[1]{{\textit{Rble}}_N(#1)}
\newcommand{\Red}[1]{{\textit{Rble}}_D(#1)}
\newcommand{\NRevLHS}[1]{\textit{Rev}_N(\mathit{LHS}(#1))}
\newcommand{\DRev}[1]{\textit{Rev}_D(#1)}
\newcommand{\DRevLHS}[1]{\textit{Rev}_D(\mathit{LHS}(#1))}
\newcommand{\WDiff}{\textit{WDiff}\xspace}
\newcommand{\RR}{\textsf{\textup{R}}\xspace}
\newcommand{\RRn}[1]{\textsf{\textup{R}{\rm(}$#1${\rm)}}}
\newcommand{\Store}{\textsf{\textup{S}}\xspace}
\newcommand{\Storen}[1]{\textsf{\textup{S}$(#1)$}}
\newcommand{\Considered}{\textsf{\textup{Considered}}\xspace}
\newcommand{\This}{\textsf{\textup{This}}\xspace}

\newcommand{\New}{\textsf{\textup{New}}\xspace}
\newcommand{\Delete}{\textsf{\textup{Delete}}\xspace}

\newcommand{\Permanent}{\textsf{\textup{P}}\xspace}
\newcommand{\SetOfRules}[1]{\mathsf{Set}(#1)}
\newcommand{\LHS}[1]{\mathsf{LHS}(#1)}
\newcommand{\needed}{\textsf{needed}\xspace}
\newtheorem{theorem}{Theorem}[section]
\newtheorem{proposition}[theorem]{Proposition}
\newtheorem{lemma}[theorem]{Lemma}
\newtheorem{corollary}[theorem]{Corollary}
\theoremstyle{definition}
\newtheorem{defn}[theorem]{Definition}
\newtheorem{note}[theorem]{Note}
\newtheorem{example}[theorem]{Example}
\newtheorem{remark}[theorem]{Remark}

\newcounter{cond}

\newcommand{\thmref}[1]{Theorem~\ref{#1}}
\newcommand{\defref}[1]{Definition~\ref{#1}}
\newcommand{\figref}[1]{Figure~\ref{#1}}

\newcommand{\remref}[1]{Remark~\ref{#1}}
\newcommand{\propref}[1]{Proposition~\ref{#1}}
\newcommand{\exref}[1]{Example~\ref{#1}}
\newcommand{\uast}{^{\textstyle \ast}} %superscript asterisk of reasonable size
\newcommand{\Astar}{A\uast}
\newcommand{\naturals}{\mathbb N}
\begin{document}

\title[Automatic Groups and Knuth--Bendix]
{Automatic Groups and Knuth--Bendix with Infinitely Many Rules}

\author{D.\ B.\ A.\ Epstein}
\address{Mathematics Institute\\ University of Warwick\\ Coventry CV4 7AL\\ UK}
\email{dbae@maths.warwick.ac.uk}

\author{Paul J.\ Sanders}
\address{Mathematics Institute\\ University of Warwick\\ Coventry CV4 7AL\\ UK}
\email{pjs@maths.warwick.ac.uk}



\thanks{Funded by EPSRC grant no. GR/K 76597}

\begin{abstract}
It is shown how to use a small finite state automaton in two variables
in order to carry out part of the Knuth--Bendix process for rewriting words in
a group. The main objective is to provide a substitute for the most
space-demanding module of the existing software which attempts to find
a \textit{shortlex}-automatic structure for a group.
The two-variable automaton can be used to store an infinite set of
rules and to carry out fast reduction of arbitrary words using this
infinite set. We introduce a new operation, which we call welding, which
applies to an arbitrary finite state automaton. In our context this operation
is vital. We point out a small potential improvement in the subset algorithm
for making a non-deterministic automaton deterministic.
\end{abstract}

\maketitle

\section{Introduction}
A celebrated result of Novikov and Boone asserts that the word problem for
finitely presented groups is, in general, unsolvable. This means that
a finite presentation of a group has been written down, with the
property that there is no algorithm whose input is a word in the
generators, and whose output states whether or not the word is
trivial. So, given a presentation of a group which one is unable to
analyze, can any help at all be given by brute force methods, using a
computer?

The answer is that some help \textit{can} be given with the kind of
presentation that arises naturally in the work of many mathematicians,
even though one can formally prove that there is no procedure that
will \textit{always} help.

There are two general techniques for trying to determine, with the help
of a computer, whether two words in a group are equal or not. One is
the Todd--Coxeter coset enumeration process and the other is the
Knuth--Bendix process. Todd-Coxeter is more adapted to finite groups
which are not too large. We are mostly interested in groups which
arise in the study of low dimensional topology. In particular they are
infinite groups, and the number of words of length $n$ rises
exponentially with $n$. For this reason, Todd--Coxeter is not much use
in practice. Well before Todd--Coxeter has had time to work
out the structure of a large enough neighbourhood of the identity in the
Cayley graph to be
helpful, the computer is out of space.

On the other hand, the Knuth--Bendix process is much better adapted to this 
task, and it has been used quite extensively, particularly by Sims,
for example in connection with computer investigations into problems
related to the Burnside problem.
It has also been used to good effect by Holt and Rees in their automated
searching for isomorphisms and homomorphisms between two given
finitely presented groups (see \cite{HR})
In connection with searching for a \textit{shortlex}-automatic
structure on a group (we say what this means in
Section~\ref{Automatic groups}),
Holt was the first person to realize that the
Knuth--Bendix process might be the right direction to choose (see
\cite {EHR}).  However, Knuth--Bendix will run
for ever on even the most innocuous hyperbolic triangle groups, which are
perfectly easy to understand. Holt's successful plan was to use
Knuth--Bendix for a certain amount of time, decided heuristically, and
then to interrupt Knuth--Bendix and use axiom-checking,
a part of automatic group theory (see
\cite[Chapter 6]{wordprocessing}), to find an automatic structure on
the group.
Thus, using the concept of an automatic group as a mechanism for
bringing Knuth--Bendix to a halt has been one of the philosophical bases
for the work done at
Warwick in this field almost from the beginning.
In addition to the works already cited in this paragraph, the reader
may wish to look at \cite{HR} and \cite{DFH}.

For a \textit{shortlex}-automatic group, a minimal set of
Knuth--Bendix rules may be infinite, but it is always a regular
language, and therefore can be encoded by a finite state machine.
In this paper, we carry this philosophical approach further,
attempting to compute this finite state machine directly, and to carry
out as much of the Knuth--Bendix process as possible using only
approximations to this machine.

Thus, we describe a setup that can handle an infinite regular set of
Knuth--Bendix rewrite rules. For our setup to be effective, we need to
make several assumptions. Most important is the assumption that we
are dealing with a group, rather than with a monoid. Secondly, our
procedures are perhaps unlikely to be of much help unless the group actually
is \textit{shortlex}-automatic.

As a computer science byproduct of our work, we produce a new
operation on automata, which we call \textit{welding}. Although this
is an operation which makes sense on the level of abstract
languages, we do not see any use for it apart from those indicated in
this paper, which is concerned very much with equations in groups.
Another computer science byproduct is a small improvement which
one can sometimes make in the process of determinizing a finite state
automaton. Since determinization is potentially exponential, even a
small improvement can be useful.

Previous computer implementations of the semi-decision procedure to
find the \textit{shortlex}-automatic structure on a group are
essentially specializations of the Knuth--Bendix procedure \cite{KB} to a
string rewriting context together with fast, but space-consuming,
automaton-based methods of performing word reduction relative to a finite set
of \textit{shortlex}-reducing rewrite rules. Since \textit{shortlex}-automaticity
of a given finite presentation is, in general, undecidable, space-efficient
approaches to the Knuth--Bendix procedure are desirable. In this paper we
present a new algorithm which performs a Knuth--Bendix type procedure relative
to a possibly infinite regular set of \textit{shortlex}-reducing rewrite rules,
together with a companion word reduction algorithm which has been designed with
space considerations in mind.

We would like to thank Derek Holt for many conversations about this
project, both in general and in detail. His help has, as always, been generous
and useful.

\section{String rewriting}
\label{strrewriting}
In this section we review some standard material on string-rewriting,
with the object of making this paper reasonably self-contained. Later
sections will review standard material on automata and automatic
groups.

\begin{defn}\label{Astar}
Let $A$ be a finite set (usually called the \textit{alphabet}).
We define $\Astar$ to be the set of strings of symbols in $A$.
In other words, $\Astar$ is the free monoid generated by $A$, with
multiplication defined by concatenation. The identity element in this
monoid is the empty string, denoted by $\epsilon$.
\end{defn}

\begin{defn}
\label{rewrite rules}
Given a finite alphabet $A$,
a subset $R$ of $A\uast \times A\uast$ is known as a
{\it rewrite system} on $A\uast$. The elements of $R$ are known as
{\it rewrite rules}.
\end{defn}

\begin{defn}
\label{elementary reduction}
\label{Thue congruence}
Associated with a rewrite system $R$ we define three relations
$\rightarrow_R$, $\rightarrow_R\uast$ and $\leftrightarrow_R\uast$ on
$A\uast$. For $u,v \in A\uast$ we write $u \rightarrow_R v$ (and say
that $u$ {\it rewrites} to $v$) if there are strings $x,y \in A\uast$
and a rewrite rule $(\lambda,\rho) \in R$ such that $u = x\lambda y$ and
$v = x\rho y$. This is also called an \textit{elementary reduction}. The
relation $\rightarrow_R\uast$ is the reflexive, transitive closure of
$\rightarrow_R$ and the relation $\leftrightarrow_R\uast$ is the
reflexive, symmetric, transitive closure of $\rightarrow_R$. The
congruence $\leftrightarrow_R\uast$ is called the {\it Thue congruence}
generated by $R$ and we denote the congruence class of a string $w \in
A\uast$ by $[w]_R$.
\end{defn}

If there is no infinite sequence $u_1\rightarrow_R u_2\rightarrow_R
\cdots$ of rewrites we say that $R$ is {\it Noetherian}. In such a
system each congruence class contains at least one {\it irreducible}
string, i.e.\ an element $w\in A\uast$ which contains no substring equal
to the left-hand side of any rewrite rule. In a Noetherian rewrite
system any string $w$ is reduced to an irreducible element of $[w]_R$ by
a finite sequence of rewrites. If each congruence class of a Noetherian
system contains a unique irreducible then the word problem in the
quotient monoid $A\uast/\leftrightarrow_R\uast$ is solved by rewriting.

A rewrite system $R$ is called
{\it confluent} if whenever $u,v,w \in A\uast$ with $u \rightarrow_R v$
and $u \rightarrow_R w$, there exists some string $z\in A\uast$ with
$v \rightarrow_R\uast z$ and $w \rightarrow_R\uast z$.
Confluence can easily be proved necessary and sufficient for each
congruence class in a Noetherian rewrite system to contain a unique
irreducible.
In this case elements of the monoid $A\uast/\leftrightarrow_R\uast$
can be defined by juxtaposition and reduction.

\subsection{Critical pair analysis.}
\label{critical case analysis}
For a finite Noetherian rewrite system $R$ the question of unique
irreducibles entails only a finite computation known as {\it critical pair
analysis}. If a finite Noetherian $R$ is not confluent one can easily prove
that the property must fail at one of a finite number of triples $(u,v,w)$.
Such triples are
obtained by considering pairs of rules in $R$ whose left-hand sides
have a non-trivial {\it overlap}. For such a pair of rules
$(\lambda_1,\rho_1)$, $(\lambda_2,\rho_2)$ there are two types of overlap.
First, a non-empty string $z$ may be a suffix of $\lambda_1=s_1z$
and a prefix of
$\lambda_2=zs_2$ (or vice versa).
Second, $\lambda_2$ may be a substring of
$\lambda_1$ (or vice versa) and we write $\lambda_1=s_1\lambda_2s_2$.

These cases are not disjoint. In particular, if one of $s_1$
and $s_2$ is trivial in the second case, it can equally well be
treated under the first case with $z$ equal either to $\lambda_1$ or
to $\lambda_2$.

\subsection{First case of critical pair analysis.}
\label{first case}
In the first case, the triple $(u,v,w)=(s_1zs_2,\rho_1s_2,s_1\rho_2)$.
There are two ways of starting to reduce $u=s_1zs_2$, namely to $v=\rho_1s_2$
and to $w=s_1\rho_2$.
Further reduction to irreducibles either gives the same irreducible
for each of the two computations, or else
gives us distinct irreducibles $v^\prime$ and $w^\prime$.
In the latter situation we can augment $R$ either with the
rule $(v^\prime,w^\prime)$ or with $(w^\prime,v^\prime)$ provided the system
obtained remains Noetherian. (If it doesn't remain Noetherian with
either choice, we will almost certainly have
to give up on the whole process.) Since $v^\prime$ was previously
congruent to $w^\prime$, the congruence on $\Astar$ is
unaltered by introducing such a rule.

Note that it is important to allow $(\lambda_1,\rho_1) = (\lambda_2,\rho_2)$
in the first case, provided there is a $z$ which is both a proper suffix and
a proper prefix of $\lambda_1=\lambda_2$.

\subsection{Second case of critical pair analysis.}
\label{second case}
In the second case, the triple $(u,v,w) = (\lambda_1,\rho_1,s_1\rho_2s_2)$.
If $\rho_1$ and $s_1\rho_2s_2$ do not reduce to the same
irreducible, we augment $R$
with a new rule $(v^\prime,w^\prime)$ or with $(w^\prime,v^\prime)$,
provided the system remains Noetherian.

\subsection{Omitting rules.}
\label{omitting rules}
In practice, it is important to remove rules which are redundant, as
well as to add rules which are essential. Omitting rules is
unnecessary in theory, provided that we have unlimited time and space at our
disposal. In practice, if we don't omit rules, we are liable to be
overwhelmed by unnecessary computation. Moreover, nearly all programs
in computational group theory suffer from excessive demands for space.
Indeed this is one of the reasons for developing the algorithms and
programs discussed in this paper. So it is important to throw away
information that is not needed and doesn't help.

For this reason, in Knuth--Bendix programs one normally looks from
time to time at each rule $(\lambda,\rho)$ to see if it can and should be
omitted. If a proper substring of
the left-hand side can be reduced, then we
are in the situation of \ref{second case}. If the two reductions
mentioned in \ref{second case}
lead to the same irreducible, we omit $(\lambda,\rho)$ from the
set of rules. If the two reductions lead to different irreducibles,
then we augment the set of rules as described in \ref{second case} and
again omit $(\lambda,\rho)$.

We also investigate whether the right-hand
side $\rho$ of a rule $(\lambda,\rho)$ is reducible to $\rho'$.
If so, we can omit $(\lambda,\rho)$ from $R$ and replace
it with the rule $(\lambda,\rho')$ without changing the congruence on $\Astar$.

\subsection{Maintaining the congruence.}
The insertion of a rule $(\lambda,\rho)$, for which we already know
that $\lambda$ and $\rho$ are congruent, clearly does not change the
congruence.

Suppose, on the other hand, that $(\lambda,\rho)$ is a rule where
$\lambda$ is known to be congruent
to $\mu$ using only rules other than $(\lambda,\rho)$. Then $(\lambda,\rho)$
can be omitted without changing the congruence.

The process of analysing critical pairs and augmenting or diminishing
the rule set without changing the congruence on $\Astar$ is known as
the {\it Knuth--Bendix Process}.
If this terminates, it gives a
finite confluent rewriting system for the congruence. Usually it
does not terminate and it produces new rules {\it ad infinitum}. At each
augmentation we have to choose one of two new
rules to insert and we have to ensure
that the augmented system is still Noetherian in order for the procedure to
continue. The choice is generally made using some total ordering on
the elements of $\Astar$.

\subsection{Ordering}
\label{ordering details}

Given an alphabet $A$, a {\it reduction ordering} on $A\uast$ is a well-ordering
which is invariant under left and right multiplication.
Suppose we have a rewriting system $R$ on $A\uast$, such that,
for each $(\lambda,\rho) \in R$,
we have $\rho < \lambda$.
Then we say that $R$ is {\it consistent with the reduction ordering}.

If $R$ is consistent with a reduction ordering, then it is Noetherian.
Moreover, to augment the rule set during the Knuth--Bendix
Completion Procedure without destroying the
Noetherian property,
we have to choose either $(v^\prime,w^\prime)$, or $(w^\prime,v^\prime)$,
using the notation introduced above where $v^\prime$ and $w^\prime$
are certain elements in $\Astar$.
We use $(v^\prime,w^\prime)$ if $v^\prime>w^\prime$ and
$(w^\prime,v^\prime)$ if $w^\prime>v^\prime$.
For more information on Knuth--Bendix and rewriting see
\cite{buchbergerloos}.

Given an ordering $<$ on a finite alphabet $A$ we can extend this to a
well-ordering on $A\uast$ by defining $ u < v$ to mean either 1) $|u| <
|v|$ or 2a) $|u| = |v|$ and 2b) $u<v$ in the
lexicographic order on $A\uast$ induced by the ordering on $A$.
This is clearly a reduction ordering on $A\uast$ and is termed the
\textit{shortlex}-ordering induced by $<$.
It is the only ordering we will use in this paper.

We order pairs $(\lambda,\rho)\in \Astar\times \Astar$ by setting
$(\lambda',\rho') > (\lambda,\rho)$ if and only if either
$\lambda'>\lambda$ or $\lambda'=\lambda$ and $\rho'>\rho$. This is
clearly a well-ordering on the set of pairs.

\subsection{Infinite sets of rules.}
In theory, critical pair analysis can be undertaken even on an
infinite set of rules $R$, provided we are working with a reduction
ordering. We set $R_1 = R$. In the $n$-th step, we
do critical pair analysis on all rules of $R_n$ such that the sum of
the lengths of left-hand side and right-hand side is at most $n$. The
effect on $R_n$ is to delete some rules, namely those such that either
the right-hand side has been shown to be reducible or the left-hand
side is reducible without using the rule itself,
and to insert others, namely those that
arise in the critical pair analysis. The
resulting set of rules is $R_{n+1}$. We can form $S = \bigcap_m
\bigcup_{n\ge m} R_n$. The congruence on
$\Astar$ induced by $R$ is the same as the congruence induced by
$R_n$ for each $n$.
It can also be proved to be the same as the congruence induced by
$S$. Moreover $S$ can be proved to be a confluent system.
The $S$-irreducibles are in
one-to-one correspondence with the elements of the monoid
$A\uast/\leftrightarrow_R\uast$.

If $R_1$ is finite, then $R_n$ is finite for each $n$.

Unfortunately, $S$ is sometimes not a
recursive set, even if the $R_n$ are all finite,
so that it cannot be computed by a Turing machine.
$S$ consists of exactly those rules with irreducible right-hand sides
and reducible left-hand sides, such that any proper substring of the
left-hand side is irreducible.

In our treatment, we will be dealing with an infinite set of rules
defined implicitly by a finite state automaton. However, we will not
attempt to perform Knuth--Bendix directly on this infinite set.

\subsection{Knuth--Bendix pass.}
\label{Knuth--Bendix pass}
One procedure for carrying out the Knuth--Bendix process is to divide the
finite set \Store of rules found so far into three disjoint subsets.
The first subset,
called {\Considered}, is the set of rules whose left-hand sides have been
compared with each other and with themselves for overlaps.
The second set of rules, called {\This}, is the set of rules waiting to
be compared with those in {\Considered}. The third set, called {\New},
consists of those rules most recently found.
Here we only sketch the process. Full details are provided in
Section~\ref{our version}.

The Knuth--Bendix process proceeds in phases, each of which is
called a \textit{Knuth--Bendix pass}. Each pass starts by looking at each
rule in {\Considered} and seeing whether it can be deleted as in
\ref{omitting rules}.
Consideration of an existing rule in {\Considered} can lead to a
new rule, in which case the new rule is added to {\New}.

Next, we look at each rule $r$ in \New to see if it is redundant.
If it is redundant it is replaced by a non-redundant rule.
The details will be given in \ref{minimizing a rule}.
The non-redundant version of the rule is moved into \This.

We then look at each rule in \This.
Its left-hand side is compared with
itself and with all the left-hand sides of rules in
{\Considered}, looking for overlaps as in \ref{first case}.
Any new rules found are added to {\New}.
Then $r$ is moved into {\Considered}.
Eventually \This becomes empty.

We then proceed to the next pass.

\section{Automata and operations on them}

This section is devoted to standard material.

\begin{defn}
A \textit{non-de\-ter\-min\-ist\-ic finite state automaton}
(abbreviated NFA)
is defined to be a quintuple
$(S, A, \mu, F, S_0)$, where $S$ is a finite set called
the \textit{set of states}, $A$ is a finite set called the
\textit{alphabet}, $\mu$ is a set of triples,
called the \textit{set of arrows},
of the form $(s,x,t)$ with $s,t \in S$ and $x \in A$ or $x=\epsilon$,
where $\epsilon$ is defined in \ref{Astar}.
$F\subset S$ is called the
\textit{set of final states} and $S_0\subset S$ is called
the \textit{set of initial states}.
The \textit{source} of an arrow $(s,x,t)$ is defined to be $s$ and the
\textit{target} is defined to be $t$.
Final states are sometimes called
\textit{accept states}, initial states are sometimes called
\textit{start states} and arrows are sometimes called
\textit{transitions}.

We define a \textit{path of arrows} in a non-de\-ter\-min\-ist\-ic automaton
$M = (S, A, \mu, F, S_0)$ to be a finite sequence of the
form $(u_0,x_1,u_1,\ldots,x_{n},u_n)$, where $n\ge 0$ and, for each
$i$ with $1\le i \le n$, $(u_{i-1},x_i,u_i) \in \mu$.
The \textit{length} of the path is $n$.
The \textit{label} associated to a path is the element $x_1\ldots x_n
\in \Astar$. If $n=0$, the label is $\epsilon \in \Astar$.
The language $L(M)$ \textit{accepted} by $M$
is defined to be the set of all labels of
paths of arrows starting with some $u_0\in S_0$ and ending with some
$u_n\in F$.
If a subset of $\Astar$ is equal to $L(M)$ for some non-de\-ter\-min\-ist\-ic
automaton $M$, then the subset is called a \textit{regular language}.
\end{defn}

\begin{defn}
A {\it partially deterministic finite state automaton} (abbreviated
PDFA) is defined to be an NFA
$M=(S,A,\mu,F,S_0)$ which contains exactly one initial state,
has no $\epsilon$-arrows and where for each $s\in S$ and $x\in A$ there is at 
most one arrow of the form $(s,x,t)$.
\end{defn}

\begin{defn}
A {\it deterministic finite state
automaton} (abbreviated DFA) is defined to be an NFA $(S,A,\mu,F,S_0)$,
in which there are no $\epsilon$-arrows, $S_0$ is a singleton
whose unique element $s_0$ is called the \textit{initial state}, and
such that, for each $s\in S$ and $x\in A$, there is exactly one arrow
of the form $(s,x,t)$.
\end{defn}

Given a non-de\-ter\-min\-ist\-ic automaton $M$ and a subset $T$ of the set of
states $S$ we define the
{\it $\epsilon$-closure} ${\mathcal E}(T)$ of $T$ to be the subset of $S$
which one can reach from some element of
$T$ by following a path of $\epsilon$-arrows.
A non-de\-ter\-min\-ist\-ic automaton can be converted into a deterministic
automaton accepting the same language as follows. The states of the
new automaton are the $\epsilon$-closed subsets of $S$ (one of these
states is the nullset). Given an $\epsilon$-closed set $T$ and $x\in
A$, we define an arrow $(T,x,P)$, where $P$ is obtained by taking the
set of targets of all $x$-arrows with source in $T$ and then taking
its $\epsilon$-closure. The initial state is the $\epsilon$-closure of
$S_0$. A state of the new deterministic automaton is final if and only
if it contains a final state of $M$.

This proves the following standard theorem.

\begin{theorem}\label{determinize}
For any NFA $M$ there is a DFA $N$ with $L(N) = L(M)$.
\end{theorem}

\begin{note}
We will use the abbreviation FSA to denote an automaton which is a
DFA or an NFA or a PDFA.
\end{note}

Computationally, the
procedure of finding the $\epsilon$-closed subsets and the arrows between
them is known as the {\it subset construction} and this is central to our word 
reduction algorithm.
There is a theoretical exponential blow-up in the subset construction
which is known to be unavoidable in general. In the cases which
come up in practice in our work, the subset construction can certainly
be a problem, but is often not as bad as the worst case analysis seems
to suggest. The implementation need only construct
those $\epsilon$-closed subsets which can be reached from $S_0$.
The space and time demands of the procedure are proportional to the
number of such subsets. We will also use lazy evaluation
to reduce the worst effects of this exponential blow-up. This
will be described later.

For a general DFA $M$, the process of finding a DFA $M^\prime$ with
$L(M^\prime) = L(M)$, such that the number of states of $M^\prime$
is minimal, is
known as {\it minimization}. The existence and uniqueness (up to isomorphism)
of such an automaton is known as the Myhill--Nerode theorem and many practical
algorithms exist to find $M^\prime$ given
$M$---for a detailed survey and comparisons see \cite{watson}.

In order to define what is meant by an automatic group we need to first 
formalize what it means for
an automaton to accept pairs of strings over an alphabet $A$.
Consider, for example, the pair of strings $(abb,ccd)$.
We regard this
pair as a string
$(a,c)(b,c)(b,d)$ over the product alphabet $A\times A$. If the pair of
strings is $(abb,ccdc)$, then we have to \textit{pad} the shorter of the two
strings to make them the same length, regarding this pair as 
the string of length four $(a,c)(b,c)(b,d)(\$,c)$.
In general, given an arbitrary pair of strings $(u,v) \in \Astar\times\Astar$,
we regard this instead as a string of pairs by adjoining a {\it padding symbol}
$\$$ to $A$ and then ``padding'' the shorter of $u$ and $v$ so that
both strings have the same length.
We obtain a string over $A\cup \{\$\} \times A\cup\{\$\}$.
The alphabet $A\cup\{\$\}$ is denoted $A^+$
and is called the {\it padded extension} of $A$.
The result of padding an arbitrary pair $(u,v)$
is denoted $(u,v)^+$. A string $w \in (A^+)\uast\times(A^+)\uast$
is called {\it padded} if there exists $u,v \in A\uast$
with $w = (u,v)^+$ (in other words, at most one of the two components
of $w$ ends with a padding symbol).

A set of pairs of strings over $A$ is called {\it regular} if the corresponding
set of padded strings is a regular language over the product
alphabet $A^+\times A^+$.

We will need two standard definitions when dealing with finite state
automata.
\begin{defn}
Let $M$ be an FSA.
We define its {\it reversal} $\Rev M$ to be the FSA obtained
from $M$ by taking the same set of states,
interchanging the subsets of initial and final states,
and then reversing the direction of all arrows.
The reversal of a DFA is in general an NFA rather than a DFA.
\end{defn}
\begin{defn}
An FSA is called {\it trim} if each state has an accepted path of
arrows passing
through it.
\end{defn}

\section{A modified determinization algorithm}
\label{A modified determinization algorithm}

In this section we discuss a modification to the usual determinization
algorithm for turning an NFA into a DFA.
Let $N$ be an NFA. The proof that $N$ can be determinized is discussed just
before the statement of Theorem \ref{determinize}.
Let $M$ be the corresponding
determinized automaton, so that
a state of $M$ is a subset of states of $N$. In practice, to find $M$,
we start with the $\epsilon$-closure of the set of initial states of
$N$ and proceed inductively. If we have found a state $s$ of $M$ as a
subset of the set of states of $N$, we fix some $x\in A$, and
apply $x$ in all possible ways to all $t\in s$, where $t$ is a state
of $N$. We then follow with $\epsilon$-arrows to form an $\epsilon$-closed
subset of states of $N$. This gives us the result of applying $x$ to
$s$. The modification we wish to make to the usual subset
construction is now explained and justified.

We will denote by $M'$ the modified version of $M$ thus obtained. $M'$ is
a DFA which accepts the same language as $M$ and $N$, but the structure of
$M'$ might be a little simpler than that of $M$.

Suppose $p$ is a state of the NFA $N$. Let $N_p$ be the same automaton
as $N$, except that the only initial state is $p$. Suppose $p$ and $q$
are distinct states of $N$ and that $L(N_p)\subset L(N_q)$.  Suppose
also that the $\epsilon$-closure of $q$ does not include $p$.  Under
these circumstances, we can modify the subset construction as follows.
As before, we start with the $\epsilon$-closure of the set of initial
states of $N$. We follow the same procedure for defining the arrows and
states of $M'$ as for $M$, except that, whenever we construct a subset
containing both $p$ and $q$, we change the subset by omitting $p$.

The situation can be generalized. Suppose that, for $1\le i \le k$,
$p_i$ and $q_i$ are states of $N$. We assume that all $2k$ states are
distinct from each other and that, for each pair $(i,j)$, the
$\epsilon$-closure of $q_i$ does not include $p_j$. Suppose further that,
for each $i$,
$L(N_{p_i})\subset L(N_{q_i})$.  We follow the same procedure for
defining the arrows and states of $M'$ as for $M$, except that, whenever
we construct a subset containing both $p_i$ and $q_i$, we change the
subset by omitting $p_i$.

\begin{theorem}\label{modified determinization}
 Under the above hypotheses, $L(M') = L(M)$.
\end{theorem}
\begin{proof}
Consider a string $w=x_1\cdots x_n\in \Astar$ which is accepted by $N$
via the path of arrows in $N$
$$(v_0,\epsilon\ast,u_1,x_1,v_1,\cdots,v_{n-1},\epsilon\ast,u_n,x_n,v_n,
\epsilon\ast,u_{n+1}).$$
This means that,
for each $i$ with $0\le i \le n$,
there is an $x_i$-arrow in $N$ from
$u_i$ to $v_i$ and
$u_{i+1}$ is in the $\epsilon$-closure of $v_i$.
Moreover $v_0$ is an initial state and $u_{n+1}$ is a final state.

Suppose inductively that after reading $x_1\cdots x_{i-1}$,
$M'$ is in state $s_{i-1}$. We assume inductively that we have
a path of arrows in $N$
$$(u_i^i,x_i,v_i^i,\epsilon\ast,u_{i+1}^i,\cdots,
u_n^i,x_n,v_n^i,\epsilon\ast,u_{n+1}^i),$$
such that $u_i^i\in s_{i-1}$ and $u_{n+1}^i$ is a final state.

The induction starts with $i=1$ and $s_0$ the initial state of $M'$.
We form $s_0$ by taking all initial states of $N$, and taking their
$\epsilon$-closure. If this subset of states of $N$ contains both
$p_j$ and $q_j$, then $p_j$ is omitted from $s_0$, the initial state
of $M'$.

If $u_1\notin s_0$, then we must have $u_1 = p_j$ for some $j$,
with $q_j\in s_0$.
Now $w\in L(N_{p_j}) \subset L(N_{q_j})$.
It follows that
we can take $u_1^1$ in the $\epsilon$-closure of $q_j$ and then define the
rest of the path of arrows for the case $i=1$. Since $q_j\in s_0$ and
$u_1^1$ is in the $\epsilon$-closure of $q_j$, $u_1^1$ is not equal to
any of the $p_r$. So $u_1^1\in s_0$ and the induction can start.

Now suppose we have a path of arrows
$$(u_i^i,x_i,v_i^i,\epsilon\ast,u_{i+1}^i,\cdots,
u_n^i,x_n,v_n^i,\epsilon\ast,u_{n+1}^i),$$
in $N$ such that $u_i^i\in s_{i-1}$ and $u_{n+1}^i$ is a final state of $N$.
We define $s_i$ from $s_{i-1}$ in the usual way,
applying $x_i$ in all possible ways to all states in $s_{i-1}$,
obtaining in particular $v_i^i$, and then taking the $\epsilon$-closure,
obtaining in particular $u_{i+1}^i$.
Finally, if, for some $r$,
$s_i$ contains both $p_r$ and $q_r$, then $p_r$ is deleted
from $s_i$ before it becomes a state of $M'$.

It now follows that either $u_{i+1}^i\in s_i$, or else, for some $r$ with
$1\le r \le k$, $u_{i+1}^i = p_r$, $q_r\in s_i$ and $p_r\notin s_i$.
In the first case we define $u_j^{i+1} = u_j^i$ and $v_j^{i+1} = v_j^i$ for
$j> i$ and the induction step is complete.
In the second case, using the fact that $x_{i+1}\cdots x_n \in L(N_{p_r})
\subset L(N_{q_r})$, we see that
we can take $u_{i+1}^{i+1}$ in the $\epsilon$-closure of $q_r$
and then define the rest of the path of arrows.
Since $q_r\in s_i$ and $u_{i+1}^{i+1}$ is in the $\epsilon$-closure of $q_r$,
$u_{i+1}^{i+1}$ is not equal to any other $p_s$
and so $ u_{i+1}^{i+1}\in s_i$. This completes the induction step.

At the end of the induction, $M'$ has read all of $w$ and is in state $s_n$.
We also have the final state $u_{n+1}^{n+1} \in s_n$, so that $w$ is accepted
by $M'$.

Conversely, suppose $w$ is accepted by $M'$. It follows easily by induction
that if $M'$ is in state $s_i$ after reading the prefix $x_1\cdots x_i$ of
$w$, then each state $u\in s_i$ can be reached from some initial state of
$N$ by a sequence of arrows labelled successively
$x_1,\ldots,x_i$, possibly interspersed with $\epsilon$-arrows.
Now $s_n$ must contain a final state, and so $w$ is accepted by $N$.
\end{proof}

\begin{remark}
The practical usage of this theorem clearly depends on having an
efficient way of determining when the condition $L(N_p) \subset L(N_q)$
is satisified. Later we will see examples of such tests which cost
virtually nothing to implement but have the potential to save an appreciable 
amount of both space and time.
\end{remark}

\section{Automatic groups}
\label{Automatic groups}
\begin{defn}
\label{automatic}
A group $G$ is called {\it automatic} if there exists a finite inverse closed
set $A$ of monoid generators of $G$ and a regular language $L$ over A
satisfying the following two properties.
\begin{enumerate}
\item The natural monoid epimorphism $\gamma : A\uast \rightarrow G$ remains
surjective when restricted to $L$.
\item The set
\begin{equation}
\label{acceptedpairs}
\{(u,v) : u,v \in L\mbox{ and }(ux)\gamma = v\gamma\mbox{ for some }
x\in A\cup\{\epsilon\}\}
\end{equation}
is regular.
\end{enumerate}
An FSA $W$ with $L(W) = L$ is called a {\it word acceptor} for $G$.
A word acceptor together
with an FSA accepting the language \eqref{acceptedpairs} is called an
{\it automatic structure} for $G$ relative to $A$.
\end{defn}
\noindent
This definition is succinct but suppresses the geometry which lies
behind the importance of this class of groups. Given a group $G$ generated by a
finite subset $A$, the {\it Cayley graph} ${\mathcal C}(G,A)$ is the graph
whose
vertices are the elements of $G$ and where an edge joins two vertices $g,h$ if
and only if there is a generator $a \in A$ with $ga = h$. Denoting images under
the natural epimorphism $A\uast \rightarrow G$ by overscores and the length of a
string $w\in A\uast$ by $|w|$,  we define a metric on ${\mathcal C}(G,A)$ by
letting
$$d(g,h) = \min\{|w| : w\in A\uast \mbox{ with }\overline w=
g^{-1}h\}.$$
\noindent
This is termed the {\it word metric}. We get the same metric by taking
each edge of the Cayley graph and giving it length one.
This makes the Cayley graph into a
geodesic space. For $w\in A\uast$ and $i \in \naturals$,
we denote by $w(i)$ the prefix
of $w$ of length $i$ (for $i \geq |w|$ this equals $w$).

Given two words $u,v \in A\uast$ and a positive real number $k$, we say
that $u$ and $v$ {\it fellow-travel} with constant $k$ in ${\mathcal C}(G,A)$
if the group elements
\begin{equation}
\label{worddiffs}
WD(G,A,u,v) = \left\{\overline{u(i)}^{-1}\overline{v(i)} : i \in \naturals
\right\}
\end{equation}
lie in the ball of radius $k$ around the identity element of the group. We then
have the following geometrical characterization of an automatic group.

% Want the next line to be \begin{theorem}[\cite[Theorem 2.3.5]{wordprocessing}]

\begin{theorem}[\cite{wordprocessing}(Theorem 2.3.5)]
Let $G$ be a group generated by a finite inverse closed set of monoid
generators $A$, and $L$ a regular language over $A$ mapping onto $G$ under the
restriction of the natural epimorphism $A\uast \rightarrow G$.
Then $L$ satisfies
property 2 of \defref{automatic} if and only if
there exists a constant $k > 0$ such that for any $u,v\in L$, if
$d(\overline u, \overline v) \leq 1$ in ${\mathcal C}(G,A)$ then $u$ and $v$
fellow-travel with constant $k$ in ${\mathcal C}(G,A)$.
\end{theorem}

\noindent
It follows immediately that, for an automatic group, the union of the sets
$WD(G,A,u,v)$
taken over all pairs $(u,v)$ with $u,v\in L$ and $d(\overline u,\overline v)\le 1$,
is a {\it finite} set. Here $d$ is the distance in the Cayley graph. This is
also the minimal length of $\overline{u^{-1}v}$ as a word over $A$.

The union of this finite
set with the set $A$ of generators plus the identity of $G$ is called the set
of {\it word differences} $WD(G,A)$ of the automatic
structure.
$WD(G,A)$ can be regarded as a PDFA over the alphabet $A^+ \times A^+$ where
an arrow labelled $(x,y) \in A^+ \times A^+$ goes from the word difference
$w_1$ to the word difference $w_2$ if and only if
$\overline x ^{-1}w_1\overline y = w_2$.
We extend the domain of
the epimorphism $A\uast \rightarrow G$ to $(A^+)\uast$
by sending the padding symbol
to the identity of $G$. The initial state is the identity of $G$ and
this is also the only final state. The resulting automaton is known as the
{\it word difference automaton} of the automatic structure.
The goal of our
main algorithm is to calculate this automaton starting from a finite
presentation of a \textit{shortlex}-automatic group (defined below). For definitive
information on automatic groups see the book \cite{wordprocessing}.

We can now
give the definition of the class of groups we are chiefly interested in.

\begin{defn}
A group $G$ is called {\it shortlex-automatic} with respect to a finite inverse
closed well-ordered set of monoid generators $(A,<)$ if
\begin{enumerate}
\item $G$ is automatic with respect to $A$.
\item A string $w\in A\uast$ is accepted by the word acceptor
if and only if $w$ is
the least element under the induced \textit{shortlex}-ordering
of $\{v : v \in A\uast\mbox{ and } \overline v = \overline w\}$ .
\end{enumerate}
\end{defn}

\section{Welding}
\label{Welding}

In this section we describe an operation, which we call
\textit{welding}, on FSAs which is central to our
Knuth--Bendix procedure. The motivation for this operation is
postponed to Section~\ref{motivation}.

\begin{defn}
\label{welding}
An FSA is called {\it welded} if it is partially deterministic, trim and has a
(partially) deterministic reversal.
These conditions imply that, given $x\in A$ and a state $t$, there is at
most one $x$-arrow with target $t$ and also that there is exactly one
initial state and one final state.
\end{defn}

Given a trim non-empty NFA $M$, we can form a welded automaton from it as 
follows.
Given any $\epsilon$-arrow $(s,\epsilon,t)$, we may identify $s$ with
$t$.
Given distinct initial states $s_1$ and $s_2$, we may identify $s_1$
with $s_2$.
Given distinct final states $t_1$ and $t_2$, we may identify $t_1$
with $t_2$.
Given distinct arrows $(s,x,t_1)$ and $(s,x,t_2)$, we may identify
$t_1$ with $t_2$.
Given distinct arrows $(s_1,x,t)$
and $(s_2,x,t)$, we may identify $s_1$ with $s_2$.
Immediately after any identification of two states, we change the set of
arrows accordingly, omitting any $\epsilon$-arrow from a state to itself.
Since the number of states continually decreases,
this process must come to an end, and at this point the automaton is
welded.

\begin{theorem}
Given a trim non-empty NFA $M$, all welded automata obtained from it by a
process like that described above are isomorphic to each other; that
is the welded automaton $Q$ is independent of the order in which
identifications are made.
Moreover $Q$ depends only on the language $L(M)$.
$Q$ is the minimal PDFA accepting $L(Q)$.
It follows that
welding can be regarded as an operation on regular languages.
\end{theorem}

\begin{proof}
For each $x\in A$, let $x^{-1}$ be its formal inverse and let $A^{-1}$
be the set of these formal inverses. We form from
$M$ an automaton over $A\cup A^{-1}$ by adjoining an arrow of the form
$(t,x^{-1},s)$ for each arrow $(s,x,t)$ of $M$, and adjoining an arrow
$(t,\epsilon,s)$ for each arrow $(s,\epsilon,t)$. We also adjoin
$(s_1,\epsilon,s_2)$ if $s_1$ and $s_2$ are either both initial states
or both final states.  We denote this new automaton by $N$.

Let $F$ be the free group generated by $A$.
We define a relation on the set of states by $s \sim t$ if
there is a path of arrows from $s$ to $t$ in $N$ whose label gives the
identity element of $F$. This is clearly an equivalence relation.
Let $Q$ be the quotient automaton, each of whose states is one of the
equivalence classes above, with arrows inherited from $M$, not
from $N$. All $\epsilon$-arrows are omitted from $Q$.
It is easy to see that $Q$ is welded.

If $M$ starts out by being welded, then it is easy to see that $Q=M$.

Consider the identifications of states
made during welding (see the passage following \defref{welding}).
It is easy to see that the equivalence classes of states used in the
definition of $Q$ are unaltered by one of these identifications.
It follows that the automaton $Q$ remains unaltered during the entire
welding process.
When no more identifications can be made, we have $Q$ itself.
This shows that $Q$ is independent of the order in which the
identifications are carried out. In fact $Q$ can be characterized as
the largest welded quotient of $M$.

We claim that every element of $L(Q)$ arises as follows. Let $(w_1,
w_2, \ldots, w_{2k+1})$ be a $2k+1$-tuple of elements of $L(M)$, where
$k\ge 0$. Now consider $w_1 w_2^{-1} \ldots w_{2k}^{-1} w_{2k+1} \in
F$, and write it in reduced form, that is, cancel adjacent formal
inverse letters wherever possible. If the result is in $A\uast$, that
is, if after cancellation there are no inverse elements,
then it is in $L(Q)$.
Moreover, any element of $A\uast$ obtained in this way is in $L(Q)$.
This is straightforward to prove. We leave the details to the reader
because we do not need the result. The proof uses the fact that $M$ is
trim.

A welded automaton is minimal.
For let $s$ and $t$ be distinct states, and let $u$ and $v$ be strings
over $A$ which lead from $s$ and $t$ respectively to the unique final
state. Then $u$ does not lead from $t$ to the final state and $v$ does
not lead from $s$ to the final state (otherwise $s$ and $t$ would be
equal). It follows that $s$ and $t$ remain distinct in the minimized
automaton.
\end{proof}

If $M$ is a non-empty trim NFA, we denote by $\Weld{M}$ the PDFA
obtained from it by welding.
To compute $\Weld{M}$ efficiently, we first
add ``backward arrows'' to $M$.
That is, for each arrow $(s,x,t)$ in $M$, including $\epsilon$-arrows,
we add the arrow $(t,x',s)$, where $x'$ represents a backwards version
of $x$. We also add $\epsilon$-arrows to connect the initial states,
and $\epsilon$-arrows to connect the final states.
We then make use of a slightly modified version of
the coincidence procedure of Sims given in
\cite[4.6]{sims}. When this stops we have a welded automaton.

In practice, in the automata which we want to weld, backward arrows are
needed in any case for some of our algorithms. The procedure described
in the preceding paragraph therefore fits our needs particularly well.

\section{A motivating example of welding}
\label{motivation}

We will look at some particular examples to see what can happen during
the Knuth--Bendix process on words in a group, and these examples will,
we hope, convince the reader of the significance of welding as introduced
in the previous section.

We will use the standard generators $x$, $y$,
and their inverses $X$ and $Y$ for the free abelian group on two
generators. Using different orderings on this set of four generators,
we will see how welding works and why we want to use it.

Consider the alphabet $A=\{x,X,y,Y\}$ with the ordering $ x < X < y < Y$,
and denote the identity of $A\uast$ by $\epsilon$.
Let $R$ be the rewriting system on $A\uast$
defined by the set of rules
$$\{(xX,\epsilon),\ (Xx,\epsilon),\ (yY,\epsilon),\ (Yy,\epsilon),\ (yx,xy),\ (yX,Xy),\ (Yx,xY),\ (YX,XY)\}.$$
Using the \textit{shortlex}-ordering on $A\uast$, it is straightforward to see
that $R$ is a confluent system.

We now change the ordering of the set of generators to
$x < y < X < Y$ and interchange the sides of the sixth rule (to get an
order reducing and therefore Noetherian system).
Once again the rules define the free abelian group on two generators.
But this time there can be no finite confluent set of rewrite
rules defining the same congruence.
To see this, we consider the set of strings
$\{xy^nX : n \in \naturals\}$. None of these is \textit{shortlex}-least within
its $\leftrightarrow_R\uast$-equivalence class.
Therefore each of these strings is reducible relative
to any confluent set of rewrite rules which defines
$\leftrightarrow_R\uast$.
On the other hand,
each proper substring of one of
the strings $xy^nX$ is clearly \textit{shortlex}-least within
its $\leftrightarrow_R\uast$-equivalence class, and is therefore
irreducible.  It follows that a confluent
set of rewrite rules must contain each of the strings $xy^nX$
as a left-hand side.
Hence, in this situation, the Knuth--Bendix procedure will never terminate.

We will show how to generate, after only a few steps,
the automaton giving the required
infinite confluent set of rewrite rules.

\begin{figure}[h]
\begin{picture}(240,50)(20,0)
\put(-12.5,20){\vector(1,0){10}}
\put(0,20){\circle{5}}
\put(-2,5){$1$}
\put(2.5,20){\vector(1,0){40}}
\put(10,30){$(x,y)$}
\put(45,20){\circle{5}}
\put(42.5,5){$2$}
\put(47.5,20){\vector(1,0){40}}
\put(55,30){$(y,y)$}
\put(90,20){\circle{5}}
\put(87.5,5){$3$}
\put(110,20){$\ldots$}
\put(140,20){\circle{5}}
\put(137.5,5){$n$}
\put(142.5,20){\vector(1,0){40}}
\put(150,30){$(y,y)$}
\put(185,20){\circle{5}}
\put(173.5,5){$n+1$}
\put(187.5,20){\vector(1,0){40}}
\put(195,30){$(y,\$)$}
\put(230,20){\circle{5}}
\put(218.5,5){$n+2$}
\put(232.5,20){\vector(1,0){40}}
\put(275,20){\circle*{5}}
\put(263.5,5){$n+3$}
\put(240,30){$(X,\$)$}
\end{picture}
\caption{The PDFA $M(r_n)$ for $n>2$.}
\label{ruletopdfa}
\end{figure}
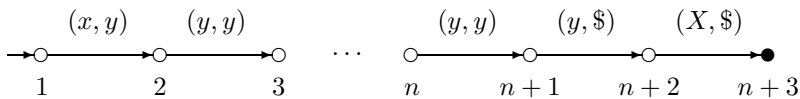

We consider the rule $r_n = (xy^nX,y^n)$
for some $n\in \naturals$. The corresponding padded string $r_n^+$ gives
rise to an $(n+3)$-state PDFA $M(r_n)$ whose accepted language consists solely
of the rule $r_n$. For $n>2$ this PDFA is shown in \figref{ruletopdfa}.

\begin{figure}
\begin{picture}(220,60)
\put(0,17.5){\circle{5}}
\put(-2.5,2.5){$1$}
\put(-12.7,17.5){\vector(1,0){10}}
\put(2.4,17.5){\vector(1,0){70}}
\put(74.8,17.5){\circle{5}}
\put(72.3,2.3){$2$}
\put(74.8,30){\circle{20}}
\put(71,37){$<$}
\put(64,46){$(y,y)$}
\put(22,22){$(x,y)$}
\put(77.2,17.5){\vector(1,0){70}}
\put(149.6,17.5){\circle{5}}
\put(147.1,2.5){$3$}
\put(105,22){$(y,\$)$}
\put(152.1,17.5){\vector(1,0){70}}
\put(224.6,17.5){\circle*{5}}
\put(221.1,2.5){$4$}
\put(179,22){$(X,\$)$}
\end{picture}
\caption{A PDFA isomorphic to $\Weld{M_n}, n > 1$.}
\label{weldrules}
\end{figure}
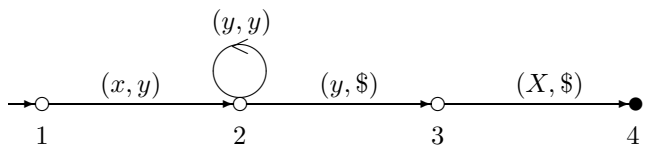

Continuing the discussion of the rules for a free abelian group on two
generators, we define $M_n$ to be the disjoint union
${\mathcal U}(M(r_1),\ldots,M(r_n))$ of the
automata $M(r_1),\ldots,M(r_n)$, with set of initial (final) states
equal to the collection of initial (final) states for the various $M(r_i)$.
If $n>1$ then $\Weld{M_n}$ is
isomorphic to the PDFA given in \figref{weldrules}, and the accepted language
of this PDFA is the set of rules $\{r_i : i \in \naturals\}$. This is
independent of $n$ if $n>1$.

So in this
example, after only two steps,
the welding procedure provides us with a PDFA whose accepted
language consists of an infinite set of rules,
each of which is a valid identity
in the group $A\uast/ \leftrightarrow_R\uast$.
Moreover, by using this PDFA to
define a suitable reduction procedure,
each of the strings $xy^nX$ with $ n\in \naturals$
can be reduced to the \textit{shortlex}-least representative of its
$\leftrightarrow_R\uast$-equivalence class.

For this group with the given ordering on the generators, it is not
hard to show that by welding the original defining rules for the group together
with the $4$ rules $\{(xyX,y),(xy^2X,y^2),(yXY,X),(yX^2Y,X^2)\}$, we obtain a
PDFA whose accepted language is a confluent set of rules.
The reduction procedure, which we will describe later, corresponding
to this PDFA will reduce {\it any} string to its
\textit{shortlex}-least representative.

\section{Rule automata}\label{Rule automata}

For the welding procedure to be used in a general Knuth--Bendix situation, we
need to show that any rules obtained are valid identities in the corresponding
monoid. We now show that if the monoid is a group (the situation
we are interested in), any rules obtained are valid identities.

\begin{defn}
\label{worddiff}
Let $A$ be a finite inverse closed set of monoid generators for a group
$G$ and, as before,  denote images under the epimorphism $(A^+)\uast
\rightarrow G$ by overscores. A {\it rule automaton for $G$} is an NFA
$M=(S,A^+~\times~A^+,\mu,F,S_0)$ together with a
function $\phi_M : S \rightarrow G$ satisfying
\begin{enumerate}
\item $F,S_0 \neq \emptyset$.
\item If $s$ is an initial or final state then $\phi_M(s)=1_G$.
\item For any $s,t \in S$ and $(x,y) \in A^+ \times A^+$ with
$(s,(x,y),t)\in\mu$
we have $\phi_M(t)~=~\overline{x}^{-1}\phi_M(s)\overline{y}$.
\item For any $s,t\in S$ with $(s,\epsilon,t)\in\mu$ we have
$\phi_M(s)=\phi_M(t)$.
\end{enumerate}
\end{defn}

Here is an equivalent way to look at the definition of a rule
automaton.
We regard $G$ as the set of
states of an automaton with alphabet $A^+~\times~A^+$ and
with an arrow $(g,(x,y),h)$ if and only if $x^{-1}gy = h$ in $G$.
Since $G$ might be infinite, this would mean considering automata with an
infinite number of states, and we would have to generalize our
definitions. (Automata with an infinite number of states are fairly
standard objects in the literature.)
In this approach, we next define what we mean by a morphism of
automata. A morphism sends states to states and arrows to arrows, but
preserves labels on arrows.
A rule automaton is then a two-variable automaton with a morphism into
the two-variable automaton $G$.
We leave the straightforward details to the reader.

\begin{example}
\label{ruleworddiff}
If $A$ is a finite inverse closed set of monoid generators for a group $G$ and
$r = (u,v)\in A\uast \times A\uast$ satisifies $\overline u = \overline v$ then,
as in \figref{ruletopdfa}, writing $r^+$ as a string
$(u_1,v_1)\cdots(u_n,v_n)\in(A^+\times A^+)\uast$, we obtain an $(n+1)$-state
rule automaton
$M(r)=(\{s_0,\ldots,s_{n}\},A^+~\times~A^+,\mu,\{s_0\},\{s_n\})$ for $G$ where
the arrows are given by
$$\mu(s_i,(u_{i+1},v_{i+1})) = s_{i+1},\, 0\leq i\leq n-1.$$
The function $\phi=\phi_{M(r)}$ assigning group elements to states is
defined inductively by $\phi(s_0) = 1_G$ and
$\phi(s_i) = \overline {u_i}^{-1} \phi(s_{i-1}) \overline{v_i}$
for $1\leq i \leq n$. As usual, the padding symbol is sent to $1_G$. The
fact that $\overline u = \overline v$ ensures that condition $2$ of
\defref{worddiff} is satisfied.
\end{example}

\begin{remark}
\label{worddiffpdfa}
For a two-variable NFA $M$ which is a rule automaton,
the PDFA $P$ obtained by applying
the subset construction to the (non-empty) set of initial states of $M$ (and
the sets that arise), is also a rule automaton for $G$ where the map $\phi_P$
is induced from $\phi_M$.
The fact that this map is well-defined follows from conditions $2,3$ and $4$
of \defref{worddiff} and the fact that $P$ is connected (by construction).

The same remark applies to the modified subset construction described
in Section~\ref{A modified determinization algorithm}.

\end{remark}
\begin{proposition}
\label{validrules}
Let $A$ be a finite inverse closed set of monoid generators for a group $G$ and
suppose that $M$ is a rule automaton for $G$. Then
\begin{enumerate}
\item Every accepted rule of $M$ is a valid identity in $G$.
\item $\Weld M$ is a rule automaton for $G$.
\end{enumerate}

Consequently every
accepted rule of $\Weld M$ is a valid identity in $G$.
\end{proposition}

\begin{proof}
Let $r = (u,v)\in A\uast\times A\uast$ be an accepted rule of $M$
and write the padded
string $(u,v)^+$ as $(u_1,v_1)\cdots(u_n,v_n)$ where
$n = \max\{|u|,|v|\}$.
Then in the PDFA $P$ obtained from $M$ (as in \remref{worddiffpdfa}), there
exists a sequence of states $s_0,\ldots,s_n$; also, for each
$i, 1\leq i\leq n$, there is a arrow from $s_{i-1}$ to $s_i$ labelled by
$(u_i,v_i)$.
Hence, from condition $3$ of \defref{worddiff}, we have
$$\phi_P(s_i) = {\overline{u_i}}^{-1}\cdots{\overline{u_1}}^{-1}
{\overline{v_1}}\cdots{\overline{v_i}}, \mbox{ for all } i \mbox{ with }
0\leq i\leq n.$$
In particular, $\overline{u_1\cdots u_n} = \overline{v_1\cdots v_n}$
and therefore the rule $r$ is valid in $G$.

To prove $2$, we need only show that when any of the operations described
just after \defref{welding} is applied to a rule automaton $M$, we
continue to have a rule automaton.
This is obvious. The final statement is now immediate.
\end{proof}

\begin{corollary}
\label{validrulescorollary}
Let $A$ be a finite inverse closed set of monoid generators for a group $G$ and
suppose that $r_1,\ldots,r_m\in A\uast \times A\uast$ are valid identities in 
$G$. Then any rule accepted by $\Weld{M(r_1),\ldots,M(r_m)}$ is also a
valid identity in $G$.
\end{corollary}

\begin{proof}
For $1\leq k\leq m$ let $M(r_k)$ be the rule automaton for $G$ as in
\exref{ruleworddiff}. Then the disjoint union
${\mathcal U}(M(r_1),\ldots,M(r_m))$
is also a rule automaton for $G$ and so the result follows by
\propref{validrules}.
\end{proof}

\begin{remark}
\label{identify}
Given a rule automaton $M$ for a group $G$, the map $\phi_M$ may not be
injective. However, if $\phi_M(s)=\phi_M(t)$ and we can somehow
determine that this is the case, then we can connect
$s$ to $t$ by an $\epsilon$-arrow, and we still have a rule automaton.
If we then weld, $s$ and $t$ will be identified. So we can hope to make
$\phi_M$ injective. However, even if $\phi_M$ is not injective, the rule
automaton $M$ can still be useful for finding equalities in the group
$G$. $M$ may not tell the whole truth, but it does tell nothing but
the truth.
\end{remark}

\section{Which words are reducible?}
\label{which words are reducible}

Suppose $G$ is a group with a finite, inverse closed and ordered set
of generators $(A,<)$.
In this section, we will work with a fixed
two-variable automaton $\Rules$.
The automaton $\Rules$ arises in our work by welding together
appropriate rules found so
far in the Knuth--Bendix process.
However we will not make use of the specific way in which $\Rules$
has been constructed.
Instead we will write down a list of properties of this automaton---when
we come to construct the automaton, it will be easy to see that the properties
are either already satisfied or that it can be arranged for them to be
satisfied.

\subsection{Properties of the rule automaton.}
\label{properties of Rules}
\begin{enumerate}
\item $\Rules$ is a trim rule automaton.
\item $\Rules$ has one initial state and one final state and they are
equal.
\item $\Rules$ and its reversal
$\Rev \Rules$ are both partially
deterministic.
These conditions imply that $\Rules$ is welded.
\item Any arrow labelled $(x,x)$, with source the initial state, also
has target the initial state.
Any arrow labelled $(x,x)$, with target the initial state, also
has source the initial state. If either of these
conditions are not fulfilled, we can
identify the source and target of the appropriate $(x,x)$-arrows,
and then weld. We will still have a rule automaton.
Later on (see Lemmas~\ref{arrows removed 1} and \ref{arrows removed 2}) we
will show that (after any necessary identifications and welding)
we can omit such arrows without loss, and, in fact, with a gain
given by improved computational efficiency.
After proving these lemmas, we will assume there are no arrows
labelled $(x,x)$ with source or target the initial state of $\Rules$.
\end{enumerate}

Since $\Rules$ is a rule automaton,
Proposition~\ref{validrules} shows that
each accepted pair $(u,v)\in L(\Rules)$ gives a valid identity $\bar u = \bar v$
in $G$.

\begin{figure}[h]
\begin{picture}(100,140)(0,-20)
\put(0,50){\circle{5}}
\put(-2,39){$1$}
\put(-12.5,50){\vector(1,0){10}}
\put(1.7,51.7){\vector(1,1){30}}
\put(-6,74){$(x,y),$}
\put(-14,64){$x>y$}
\put(1.7,48.3){\vector(1,-1){30}}
\put(-10,27){$(x,y),$}
\put(-2,17){$x<y$}
\put(33.5,83.5){\circle*{5}}
\put(31.5,73){$2$}
\put(33.5,96){\circle{20}}
\put(29.8,103){$<$}
\put(23,112){$(x,y)$}
\put(33.5,16.5){\circle{5}}
\put(31.5,21){$3$}
\put(33.5,4){\circle{20}}
\put(29.8,-8){$<$}
\put(23,-17){$(x,y)$}
\put(67,50){\circle*{5}}
\put(65,39){$4$}
\put(2.5,50){\vector(1,0){62}}
\put(22,54){$(x,\$)$}
\put(35.2,81.8){\vector(1,-1){30}}
\put(52,70){$(x,\$)$}
\put(35.2,18.2){\vector(1,1){30}}
\put(53,25){$(x,\$)$}
\put(67,50){\vector(1,0){62}}
\put(89,54){$(x,\$)$}
\put(131.5,50){\circle*{5}}
\put(129.5,39){$5$}

\end{picture}
\caption{The automaton $\SLtwo$. Solid dots represent final states.
Roman letters represent arbitrary
letters from the alphabet $A$ and the labels on the arrows indicate multiple
arrows.
For example, from state $2$ to itself there is one arrow for each pair in
$A\times A$.}

\label{SL2}
\end{figure}
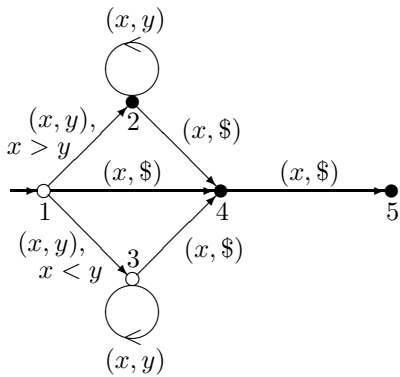

The automaton $\Rules$ may accept pairs $(u,v)$ such that $u$ is shorter
than $v$. We cannot consider such a pair as a rule and so we want to
exclude it. To this end we introduce the automaton $\SLtwo$. This is a five
state automaton, depicted in Figure~\ref{SL2}, which accepts pairs $(u,v)
\in \Astar\times\Astar$, such that $u$ and $v$ have no common prefix, $u$ is 
\textit{shortlex}-greater than $v$ and $|v|
\le |u| \le |v|+2$. By combining $\SLtwo$ with $\Rules$, we obtain a regular
set of rules
$\SetOfRules \Rules$, which is possibly infinite, namely
$L(\Rules)\cap L(\SLtwo)$. An automaton accepting this set can be constructed
as follows. Its states are pairs $(s,t)$, where $s$ is a state of $\Rules$ and
$t$ is a state of $\SLtwo$. Its unique initial state is the pair of initial
states in $\Rules$ and $\SLtwo$. A final state is any state $(s,t)$ such that
both $s$ and $t$ are final states. Its arrows are labelled by $(x,y)$, where
$x\in A$ and $y\in A^+$. Such an arrow corresponds to a pair of arrows, each
labelled with $(x,y)$, the first from $\Rules$ and the second from $\SLtwo$.

\subsection{Restrictions on relative lengths.}
\label{right-hand sides at most two shorter}
The restriction $|u| \le |v| +2$ needs some explanation. The point is that if
we have a rule with $|u| > |v| +2$, then we have an equality $\bar u = \bar v$
in $G$. We write $u = u'x$, where $x\in A$. The formal inverse $X$ of $x$
is also an element of $A$. We therefore have a pair
$(u', vX)$ which represent equal elements in $G$. If our set of rules were
to contain such a rule, then $u=u'x$ would reduce to $vXx$, and this reduces
to $v$, making the rule $(u,v)$ redundant. This leads to an obvious technique
for transforming any rule we find into a new and better
rule with $|v| \le |u| \le |v|+2$. Since we take this into account when
constructing the automaton $\Rules$, we are justified in making the
restriction.

This analysis can be carried further. Let
$u=u_1\cdots u_{r+2}=u'u_{r+2} = u_1 u''$ and
let $v=v_1\cdots v_r$. If $u_1>v_1$, then the rule $(u,v)$ can be
replaced by the better rule $(u',vu_{r+2}^{-1})$.
If $u_2 > u_1^{-1}$, then $(u,v)$ can be replaced by $(u'',u_1^{-1}v)$.
We do in fact carry out these steps when installing new rules, but we
have not so far tried adjusting the
finite state automaton $\SLtwo$ accordingly to see what effect this
would have on the whole process.

\subsection{Rules for which no prefix or suffix is a rule.}
At the moment, it is possible for an element $(u,v)^+$ of $\SetOfRules\Rules$
to have a prefix or suffix which is also a rule. This is undesirable because
it makes the computations we will have to do bigger and longer without any
compensating gain.

Recall that the automaton recognizing $\SetOfRules\Rules$
is the product of $\Rules$ with $\SLtwo$, the initial state being the product
of initial states and the set of final states being any product of final states.

We remove from $\Rules$ any arrow labelled $(x,x)$ from the initial state to
itself. We then form the product automaton, as described above, with two
restrictions. Firstly, we omit any arrow whose source is a product of final
states. Secondly, we omit completely the state and all arrows whose source or
target is the state with first component equal to $s_0$, the initial state of
$\Rules$, and second component equal to state $3$ of $\SLtwo$
(see \figref{SL2}). We call the resulting automaton $\Rules^\prime$.

\begin{lemma}\label{arrows removed 1} The language accepted by
$\Rules^\prime$ is the set of labels
of accepted paths in the product automaton, starting from the product
of initial states and ending at a product of final states, such that the only
states along the path with first component equal to $s_0$ are at the
beginning and end of the path.
\end{lemma}
\begin{proof} First consider an accepted path $\alpha$ in $\Rules'$.
The only arrows in $\Rules'$ with source having first component $s_0$
are those with source the product of initial states.
In $\SLtwo$ it is not possible to
return to the initial state. It follows that $\alpha$ has the required form.

Conversely any such path in the product automaton also lies in $\Rules'$
because it avoids all omitted arrows.
\end{proof}

\begin{lemma}\label{arrows removed 2} The language accepted by
$\Rules^\prime$ is the subset of $\SetOfRules\Rules$ which has no proper
suffix or proper prefix in $\SetOfRules\Rules$.
\end{lemma}
\begin{proof}
If $\alpha$ is an accepted path in $\Rules'$, then it is clearly in
$\SetOfRules\Rules$. Moreover if it had a proper suffix or proper prefix
which was in $\SetOfRules\Rules$, there would be a state in the middle
of $\alpha$ with first component $s_0$. We have seen that this is impossible
in Lemma~\ref{arrows removed 1}.

Conversely, we must show that if $\alpha$ is an accepted path in
the product automaton such that no proper prefix and no proper suffix of
$\alpha$ would be accepted by the product automaton, then no state
met by $\alpha$, apart from its two ends, has $s_0$ as a first component.
Let $\alpha = ((s_0,1),(u_1,v_1), q_1, \ldots, (u_n,v_n), q_n)$,

First suppose $u_1 < v_1$. Since $\alpha$ is accepted by $\SLtwo$,
we must have $v_n = \$$.
Let $r<n$ be chosen as large as possible so that the first component
of $q_r$ is $s_0$.
Then $(u_{r+1},v_{r+1})\ldots (u_n,v_n)$ will be accepted by $\Rules$
and will be accepted by $\SLtwo$ because $v_n= \$$.
Since this cannot be a proper suffix of $\alpha$ by assumption, we must have 
$r=0$. Hence $q_i$ has a first component equal to $s_0$ if and only if
$i=0 $ or $i=n$.

Next note that we cannot have $u_1=v_1$. This is because there is no
arrow labelled $(u_1,u_1)$ in $\SLtwo$ with source the initial state,
so $\alpha$ would not be accepted by the product automaton.

Now suppose that $u_1>v_1$ and let $r>0$ be chosen as small as
possible, so that the first component of $q_r$ is $s_0$.
Since $u_1>v_1$, the second component of $q_r$ will be a final state
(see \figref{SL2}). Since $\alpha$ has no accepted proper prefix, we
must have $r=n$. Hence $q_i$ has a first component equal to $s_0$ if
and only if $i=0 $ or $i=n$.

So we have proved the required result for each of the three
possibilities.
\end{proof}

Let $w=x_1\cdots x_n\in A\uast$ be a string which we wish to reduce to a
$\SetOfRules\Rules$-irreducible.
It is important for this to be done quickly, as it has been observed by many
people that the Knuth--Bendix process for strings spends most of its time
reducing. Reduction needs to be carried out during critical pair analysis.

Reduction with respect to $\SetOfRules\Rules$
is done in a number of steps. First we find the shortest reducible
prefix of $w$, if this exists. Then we find the shortest suffix of that which
is reducible. This is a left-hand side of some rule in $\SetOfRules\Rules$.
Then we find the corresponding right-hand side and substitute this for the
left-hand side which we have found in $w$. This reduces $w$ in the
\textit{shortlex}-order. We then repeat the operation until we obtain an
irreducible string. The process will be described in detail in this section
and in the subsequent two sections. There is an outline of the reduction
process in \ref{history stack}.

Our objective in this section
is to find the shortest reducible prefix of $w$, if this exists.
To achieve this, we must determine whether $w$ contains a
substring which is the left-hand side of rule belonging to $\SetOfRules\Rules$.

Let $\Rules^{\prime\prime}$ be the automaton obtained
from $\Rules^\prime$ (see Lemmas~\ref{arrows removed 1} and 
\ref{arrows removed 2}) by adding
arrows labelled $(x,x)$ from the initial state to the initial state.

We construct an NFA $\NRed\Rules$ in one variable
by replacing each label of the form $(x,y)$
on an arrow of $\Rules^{\prime\prime}$ by $x$.
Here $x\in A$ and $y\in A^+$.
The name of the automaton $\NRed\Rules$
refers to the fact that the automaton accepts reducible
strings, and does so non-deterministically.
We obtain an NFA with no $\epsilon$-arrows.
However there may be many arrows labelled $x$ with a given source.
Let $\LHS\Rules$ be the regular language of left-hand sides of rules in
$\SetOfRules\Rules$ such that no proper prefix or proper suffix of the rule
is itself a rule.

\begin{lemma}
$\Astar.\LHS\Rules= L(\NRed\Rules)$.
\end{lemma}
\begin{proof}
Because of the extra arrows labelled $(x,x)$ from initial state to
initial state, inserted into $\Rules^{\prime\prime}$, the inclusion
$\Astar.\LHS\Rules\subset L(\NRed\Rules)$ is clear.

If $u$ is accepted by $\NRed\Rules$, there is a corresponding pair
$(u,v)$ accepted by $\Rules^{\prime\prime}$. We find a maximal common prefix
$p$ of $u$ and $v$, so that $u=pu'$ and $v=pv'$.
$\Rules^{\prime\prime}$ remains in the initial state while reading $(p,p)$.
Since the initial state of $\SLtwo$ is not a final state, $(u',v')$
must be non-empty. Since there is no way of
returning to the initial state of $\SLtwo$, once $\Rules^{\prime\prime}$ starts
reading $(u',v')$, it can never return to the initial state, and therefore
$(u',v')$ must be accepted by $\Rules^\prime$. Therefore $u'\in \LHS\Rules$,
as claimed.
\end{proof}

To find the shortest reducible prefix of a given string $w$
we could feed $w$ into the FSA $\NRed\Rules$. However, reading a string
with a non-deterministic automaton is very time-consuming, as all
possible alternative paths need to be followed.

For this reason, it may at first sight seem sensible to
determinize the automaton. However, determinizing a non-deterministic
automaton potentially leads to an exponential increase in size.
The states of the determinized automaton are subsets of the non-deterministic
automaton, and there are potentially $2^n$ of them if there were $n$
states in the non-deterministic automaton.
By trying examples, we have observed that
the theoretical exponential blow-up in this construction
is sometimes a practical reality for the automaton $\NRed\Rules$.

For this reason, we use a {\it lazy
state-evaluation} form of the subset construction. The lazy evaluation strategy
(common in compiler design---see for example \cite{aho}), calculates the
arrows and subsets as and when they are needed, so that a gradually increasing
portion $P(\Rules)$ of the determinized version $\Red\Rules$ of
$\NRed\Rules$ is all that
exists at any particular time.

Lazy evaluation is not automatically an advantage. For example, if in the
end one has to construct virtually the whole determinized automaton
$\Red\Rules$ in any
case, then nothing would be lost by doing this immediately. In our special
situation,
lazy evaluation \textit{is} an advantage for two reasons. First,
during a single pass of the Knuth--Bendix process
(see Paragraph~\ref{Knuth--Bendix pass}), only a
comparatively small part of the determinized one-variable automaton
$\Red\Rules$ needs to be constructed.
In practice, this phenomenon is particularly
marked in the early stages of the computation, when the automata are
far from being the ``right'' ones. Second,
this approach gives us the
opportunity to abort a pass of Knuth--Bendix, recalculate on the basis
of what has been discovered so far in this pass, and then restart the
pass. If an abort seems advantageous early in the pass, very little
work will have been done in making the structure of the determinized
version of $\Red\Rules$ explicit.

We now describe the details of this strategy.

At the start of a Knuth--Bendix pass (see
Paragraph~\ref{Knuth--Bendix pass})
we let $P(\Rules)$ be the one-variable automaton consisting
of a single non-final state
containing only the ordered pair of initial states of $\Rules$
and $\SLtwo$. At a subsequent time during the pass, $P(\Rules)$ may have increased,
but it will always be a portion of $\Red\Rules$.

Suppose now that we wish to find the shortest prefix of the string
$w=x_1\cdots x_n\in A\uast$ which is $\SetOfRules\Rules$-reducible.
Suppose that $s_0,s_1,\ldots,s_k$ are states of $P(\Rules)$, where $0\leq k\leq n-1$,
that $s_0$ is the start state of $P(\Rules)$,
and that,
for each $i$ with $1\leq i\leq k$, we have $\mu(s_{i-1},x_i) = s_i$.
Suppose that
the target of the arrow $\mu(s_k,x_{k+1})$ is not yet defined.

By definition, the subset construction applied to the state $s_k$ of
$P(\Rules)$ under the alphabet symbol $x_{k+1}$ yields the set
$\mu_1(s_k,x_{k+1})$ as follows. For each $(s',t')\in s_k$, we look
for all arrows in $\NRed\Rules$
labelled $x_{k+1}$ with source $(s',t')$.
If $(s,t)$ is the target of such an arrow, then $(s,t)$ is an element
of $\mu_1(s_k,x_{k+1})$.
Note that this subset is always non-empty, because the initial
state of $\NRed\Rules$ is an element of each $s_i$.

In the standard determinization procedure
one would now look to see whether there is already a state $s_{k+1}$ of
$P(\Rules)$ which is equal to $\mu_1(s_k,x_{k+1})$.
If not one would create such a state $s_{k+1}$.
One would then insert an
arrow labelled $x_{i+1}$ from $s_k$ to $s_{k+1}$, if there wasn't
already such an arrow.
A new state is defined to be a final state of $P(\Rules)$ if and
only if the subset contains a final state of $\NRed\Rules$.

Of course, one does not need to determine the subset $\mu_1(s_k,x_{k+1})$
if there is already an arrow in $P(\Rules)$ labelled $x_{k+1}$ with
source $s_k$, because in that case the subset is already computed and stored.

In our procedure we improve on the procedure just
described. The point is that
$\mu_1(s_k,x_{k+1})$ may contain pairs which are
not needed and can be removed. From a practical point of view
this has the advantage of saving space and reducing the amount of computation
involved when calculating subsequent arrows. Specifically, we remove a
pair $(p,q^\prime)$ from $\mu_1(s_k,x_{k+1})$ if
$q^\prime$ is state $3$ of $\SLtwo$ (see \figref{SL2}) and
$\mu_1(s_k,x_{k+1})$ also contains the pair $(p,q)$
where $q$ is state $2$ of ${\rm SL}(5,A)$. Removing all such pairs
$(p,q^\prime)$ yields the set $\mu_P(s_k,x_{k+1})$ and we add the
corresponding arrow and state to $P(\Rules)$,
creating a new state if necessary. We make the state a final
state if the subset contains a final state of $\NRed\Rules$. The
validity of this modification follows from
Theorem~\ref{modified determinization}, and we see that
some prefix of $w$ arrives at a final state of $P(\Rules)$ if and only if
$w$ is $\SetOfRules\Rules$-reducible.

When finding the corresponding left-hand side of a rule inside $w$,
we need never compute beyond a final state of $P(\Rules)$.
As a space-saving and time-saving measure our
implementation therefore replaces each final state of $P(\Rules)$, as soon as
it is found, by the empty set of states. As remarked above, the standard
determinization
of $\NRed\Rules$ never produces an empty set of states, and so there is no
possibility of confusion.

Reading $w$ can be quite slow if many states need to be added to $P(\Rules)$
while it is being read. However, reading $w$ is fast when no states
need to be built. In practice, fairly soon after a Knuth--Bendix pass
starts, reading becomes rapid, that is, linear with a very small constant.

\section{Finding the left-hand side in a string}
\label{finding the left-hand side}
We retain the hypotheses of Section~\ref {which words are reducible}.
Namely, we have a two-variable automaton $\Rules$ satisfying the
conditions of Paragraph~\ref{properties of Rules}. We are given a word
$w=x_1\cdots x_n$, and we wish to reduce it. In the previous section
we showed how to find the minimal reducible prefix $w'=x_1\cdots x_m$
of $w$ with respect
to the rules implicitly specified by $\Rules$. We now wish to find the
minimal suffix of $w'$ which is a left-hand side of some rule in
$\SetOfRules\Rules$. The procedure is quite similar to that of the
previous section.

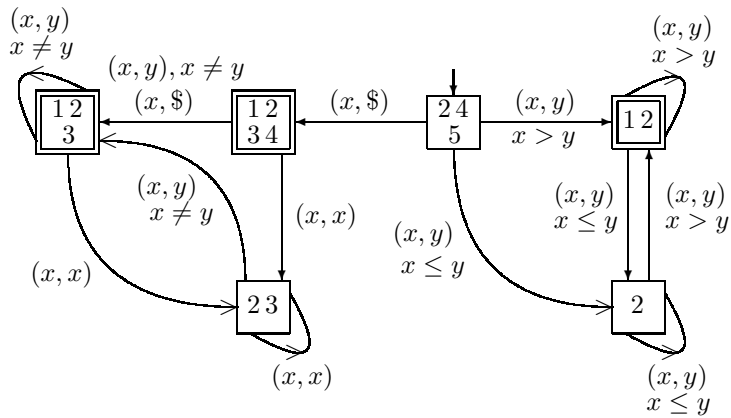
\begin{figure}[h]
\begin{picture}(250,160)(0,-50)
\put(160,80){\vector(0,-1){10}}
\put(150,50){\grid(20,20)(20,20)}
\put(154,62){$2\,4$}
\put(158,52){$5$}
\put(170,60){\vector(1,0){50}}
\put(183,65){$(x,y)$}
\put(182,52){$x>y$}
\put(220,50){\grid(20,20)(20,20)}
\put(222,52){\grid(16,16)(16,16)}
\put(224,58){$1\,2$}
\qbezier(225,70)(260,90)(240,55)
\put(239,74.5){$>$}
\put(235,93){$(x,y)$}
\put(235,83){$x>y$}
\put(220,-20){\grid(20,20)(20,20)}
\put(226,50){\vector(0,-1){50}}
\put(198,30){$(x,y)$}
\put(198,20){$x\leq y$}
\put(234,0){\vector(0,1){50}}
\put(240,30){$(x,y)$}
\put(240,20){$x>y$}
\put(227,-12){$2$}
\qbezier(225,-20)(260,-40)(240,-5)
\put(238,-29.5){$>$}
\put(233,-39){$(x,y)$}
\put(233,-49){$x\leq y$}
\qbezier(160,50)(160,-10)(220,-10)
\put(213,-12.5){$>$}
\put(137,15){$(x,y)$}
\put(140,3){$x\leq y$}

\put(150,60){\vector(-1,0){50}}
\put(76,48){\grid(24,24)(24,24)}
\put(78,50){\grid(20,20)(20,20)}
\put(82,62){$1\,2$}
\put(82,52){$3\,4$}
\put(78,-20){\grid(20,20)(20,20)}
\put(82,-12){$2\,3$}
\qbezier(83,-20)(118,-40)(98,-5)
\put(96,-29.5){$>$}
\put(91,-39){$(x,x)$}
\put(2,48){\grid(24,24)(24,24)}
\put(4,50){\grid(20,20)(20,20)}
\put(76,60){\vector(-1,0){50}}\
\put(8,62){$1\,2$}
\put(12,52){$3$}
\qbezier(2,53)(-18,92)(21,72)
\put(-2,76){$<$}
\put(-8,96){$(x,y)$}
\put(-8,86){$x\neq y$}
\qbezier(26,53)(81,53)(81,0)
\put(25.5,50.5){$<$}
\qbezier(14,48)(14,-10)(78,-10)
\put(71,-12.5){$>$}
\put(95,48){\vector(0,-1){48}}
\put(100,22){$(x,x)$}
\put(113,65){$(x,\$)$}
\put(39,65){$(x,\$)$}
\put(29,78){$(x,y),x\neq y$}
\put(40,33){$(x,y)$}
\put(45,23){$x\neq y$}
\put(0,0){$(x,x)$}
\end{picture}
\caption{This PDFA arises by applying the accessible subset construction to 
$\Rev\SLtwo$ in the case where the base alphabet has more than one element.
Each state is a subset of the state set of $\Rev\SLtwo$ and final states have 
a double border. This PDFA, when reading a
pair $(u,v)$ from right to left, keeps track of whether $u$ is longer than $v$
or not, which it discovers immediately since padding symbols if any must occur
at the right-hand end of $v$. Note that this automaton is minimized.} 
\label{RevSLsubsets}
\end{figure}

We form the two-variable
automaton $\Rev\Rules$, which we combine with $\Rev\SLtwo$.
The first automaton is, by hypothesis, partially deterministic.
If we determinize the second automaton, we obtain another PDFA.
\figref{RevSLsubsets} shows the determinization of $\Rev\SLtwo$,
where the subsets of states of $\SLtwo$ are explicitly recorded. 

We take the product of the two automata
$\Rev\Rules$ and $\Rev\SLtwo$. A new state is a pair of old states.
An arrow is a pair of arrows with the same label $(x,y)$.
The initial state in the product is the unique pair of initial states.
A final state in the product is a pair of final states.

To form the one-variable non-deterministic automaton $\NRevLHS\Rules$
without $\epsilon$-arrows, we use the same states and arrows as in the
product automaton, but replace each label
of the form $(x,y)$ in the product automaton by
the label $x$.
The deterministic one-variable automaton $\DRevLHS\Rules$
can then be constructed using the subset construction.

We have gone through the above description to give the reader a theoretical
understanding of what is going on before going into details.
Also, our procedure may be
adaptable to related situations which are not identical to this one.
In fact, we
use not the construction just described, but
a related construction which we describe below.
The point of what we do may not become fully apparent until we get to
Section~\ref{finding the right-hand side}.

\subsection{Reversing the rules.}
\label{structure} We first describe a two-variable PDFA $M$ which accepts
exactly the reverse of each rule $(\lambda,\rho)^+$ in $\SetOfRules\Rules$
such that no proper suffix and no proper prefix of $(\lambda,\rho)^+$ is in
$\SetOfRules\Rules$ (cf. Lemma~\ref{arrows removed 2}).
We assume that we have a two-variable automaton $\Rules$ satisfying the
conditions of Paragraph~\ref{properties of Rules}.

A state of $M$ is a triple $(s,i,j)$, where $s$ is a state of $\Rev\Rules$,
$i\in\{0,1,2\}$ and $j\in\{+,-\}$.
$M$ has a unique initial state $(s_0,0,+)$ where $s_0$
is the unique initial state of $\Rev\Rules$. In addition,
$M$ has three final states $f_0=(s_0,0,-), f_1=(s_0,1,-)$ and $f_2=(s_0,2,-)$.
We do not allow states of $M$ of the form $(s_0,i,j)$, except for the initial
state and the three final states just mentioned. We will construct the arrows
of $M$ to ensure that
any path of arrows accepted by $M$ has first component equal to $s_0$ for
its initial state and its final state and for no other states. (Compare this 
with Lemma~\ref{arrows removed 1}.)

The intention is that in a state $(s,i,j)$, $i$ represents the number of 
padded symbols occurring in any path of arrows from the initial 
state of $M$ to $(s,i,j)$. By \ref{right-hand sides at most
two shorter}, the padded symbols must be of the form $(x,\$)$, where $x\in A$.
Furthermore, there are zero, one or two padded symbols in any rule, and, if
padded symbols appear, they are at the right-hand end of a rule. This means
that they are the first symbols read by $M$.

The $j$ component is
intended to represent whether an arrow with source $(s,i,j)$ is
permitted with label a  padded symbol.
We take $j=+$ if a padded symbol is permitted, and $j=-$ if a padded
symbol is not permitted.

The following conditions determine the arrows in $M$.
\begin{enumerate}
\item Each arrow of $M$ is labelled with some $(x,y)$, where
$x\in A$ and $y\in A^+$.
\item $(s,i,j)^{(x,\$)}$ is defined if and only if 1) $t=s^{(x,\$)}$ is defined
in $\Rev\Rules$, and 2a) $(s,i,j)=(s_0,0,+)$, the initial state, or
2b) $(i,j)=(1,+)$. In case 2a) the
target is $(t,1,+)$,
unless $t$ is the final state of $\Rev\Rules$, in which case
the target is $f_1 = (s_0,1,-)$. In case 2b), the target is $(t,2,-)$,
which may possibly be equal to $f_2$.
The final state $f_1$ arises in case 2a) when we have a rule $(x,\epsilon)$,
which means that the generator $x$ of our group represents the trivial element.
The final state $f_2$ arises in case 2b) when we have a rule 
$(x_1x_2,\epsilon)$. This kind of rule arises when $x_1$ and $x_2$ are inverse 
to each other, usually formal inverses.
\item There are no arrows with source $f_i$.
\item Suppose $(s,i,j)$ is not a final state.
\label{non-padded M transitions}
Then $(s,i,j)^{(x,y)}$ with $x,y\in A$ is 
defined if and only if
1) $t=s^{(x,y)}$ is defined in $\Rev\Rules$, and 2) if $t=s_0$ then
2a) $i=0$ and $x>y$ or 2b) $i>0$ and $x\neq y$.
We then have $(s,i,j)^{(x,y)} = (t,i,-)$. 
This condition corresponds to the requirement that $(u,v)$ can only be a rule
if a) $u$ and $v$ have the same length and $u_1 > v_1$, where these are the
first letters of $u$ and $v$ respectively, or b)
if $u$ is longer than $v$ and $u_1\neq v_1$.
\end{enumerate}

\begin{lemma}
The language accepted by $M$ is the set of reversals of rules
$(\lambda,\rho)^+ \in \SetOfRules\Rules$ such that no proper suffix
and no proper prefix of
$(\lambda,\rho)^+$ is in $\SetOfRules\Rules$.
\end{lemma}
The proof of this lemma is much the same as the proofs of
Lemmas~\ref{arrows removed 1} and \ref{arrows removed 2}.
We therefore omit it.

Using the above description of $M$, we now describe how to obtain a
non-deterministic
one-variable automaton $\NRevLHS\Rules$ from $M$ in an analogous
manner to that
used to obtain $\NRed\Rules$ from $\Rules''$ in 
Section~\ref{which words are reducible}.
$\NRevLHS\Rules$ accepts reversed left-hand sides
of rules in $\SetOfRules\Rules$ which do not have a proper prefix or a proper
suffix which is in $\SetOfRules\Rules$.
$\NRevLHS\Rules$ has the same set of states as $M$ and the same set of
arrows. However, the label $(x,y)$ with $x\in A$ and $y\in A^+$
of an arrow in $M$ is replaced by the label $x$ in
$\NRevLHS\Rules$
The two automata, $N$ and $\NRevLHS\Rules$, have the same initial state
and the same final
states. Hence $\NRevLHS\Rules$ accepts all reversed left-hand
sides of rules whose reversals are accepted by $M$.

The one-variable automaton $Q(\Rules)$ is
formed from $\NRevLHS\Rules$ by a modified subset construction, using lazy
evaluation. $Q(\Rules)$ is part of the one-variable
PDFA $\DRevLHS\Rules$, the determinization of $\NRevLHS\Rules$. As we shall
see, a string is accepted by $Q(\Rules)$ only if its reversal $\lambda$ is the 
left-hand side of a rule in $\SetOfRules\Rules$ and no proper substring of
$\lambda$ has this property. 

\begin{note}
In order to
construct states and arrows in $Q(\Rules)$, one only needs to have 
access to $\Rev\Rules$, that is, neither $M$ nor $\NRevLHS\Rules$ has to be 
explicitly constructed. 
\end{note}

\subsection{The algorithm for finding the left-hand side.}
\label{finding lhs algorithm}
Suppose we have a string $x_1\cdots x_n \in \Astar$ and we know it has
a suffix which is the left-hand side of some rule in $\SetOfRules\Rules$.
Suppose no proper prefix of $x_1\cdots x_n$ has this property.
We give an algorithm that finds the shortest such suffix.
 
We read the string from right to left, starting with $x_n$.
We assume that $x_{k+1} x_{k+2} \cdots x_n$ has been read so far and that
as a result the current state of $Q(\Rules)$ is $S_k$, where
$S_k$ is a state of $Q(\Rules)$ (so $S_k$ is a subset of the set of states of 
$\NRevLHS\Rules$).

We start the algorithm
with $k=n$ and the current state of $Q(\Rules)$ equal
to the singleton $\{(s_0,0,+)\}$ whose only element is the initial state
of $M$,
where $s_0$ is the initial state of $\Rev\Rules$. 
$Q(\Rules)$ has three final states, namely the singleton sets
$\{f_i\}$ for $i=0,1,2$.

The steps of the algorithm are as follows:
\begin{enumerate}
\item\label{start} Record the current state as the $k$-th entry in an
array of size $n$, where $n$ is the length of the input string.
\item If the current state is not a final state, go to Step~\ref{jump}.
If the current state is a final state, then stop.
Note that the initial state of $Q(\Rules)$
is not a final state, so this step does not
apply at the beginning of the algorithm.
If the current state is a final state, then the
shortest suffix of $x_1\cdots x_n$ which is the left-hand side of a rule
in $\SetOfRules\Rules$ can then be proved to be $x_{k+1} x_{k+2} \cdots x_n$. 
\item\label{jump}
If the arrow labelled $x_k$ with source the current state is
already defined, then redefine the current state to be the target of this arrow
and decrease $k$ by one.
\item If the preceding step does not apply, we have to compute the target
$T$ of the arrow labelled $x_k$ with source the current state $S_k$.
We do this by looking for all arrows labelled $x_k$
in $\NRevLHS\Rules$ with source in $S_k$, and define $T$ to be the set of
all targets. Note that this set of targets cannot be empty since we know
that some suffix of $x_1\cdots x_n$ is accepted by $\NRevLHS\Rules$.
\item\label{modified} There are two modifications which we can make to the 
previous step. 
\begin{enumerate}
\item\label{modifieda} Firstly, if the set of targets contains some final state
$f_j$, then we look for the largest value of $i=0,1,2$ such that $f_i\in T$
and redefine $T$ to be $\{f_i\}$.
We then insert into $Q(\Rules)$
an arrow labelled $x_k$ from $S_k$ to this final state.
If we have found that $T$ is a final state, we set $S_{k-1}$ equal to $T$, 
decrease $k$ by one, and go to Step~\ref{start}.
\item\label{modifiedb}
Secondly, if, while calculating the set $T$, we find that a state $s$ of
$\Rev\Rules$ occurs in more than one triple $(s,i,j)$, then we only include
the triple with the largest value of $i$. For this to be well-defined,
we need to
know that $(s,i,+)$ and $(s,i,-)$ cannot both come up as potential
elements of
$T$---this is addressed in the proof of Theorem~\ref{lhs algorithm} along with 
justifications of the other modifications.
\end{enumerate}
\item Having found $T$, see if it is equal to some state $T'$ of $Q(\Rules)$
which has already been constructed. If so, define an arrow labelled $x_k$
from $S$ to $T'$.
\item If $T$ has not already been constructed,
define a new state of $Q(\Rules)$
equal to $T$ and define an arrow labelled $x_k$ from $S$ to $T$.
\item 
Set the current state equal to $T$ and decrease $k$ by one.
Then go to Step~\ref{start}.
\end{enumerate}

\begin{theorem}\label{lhs algorithm} Suppose $x_1\cdots x_n$ has a suffix
which is the left-hand side of a
rule in $\SetOfRules\Rules$ and suppose no prefix of $x_1\cdots x_n$
has this property.
Then the above algorithm correctly computes the shortest such
suffix.
\end{theorem}
\begin{proof}
We first show that the modification in Step~\ref{modifiedb} is well-defined in
the sense that triples $(s,i,+)$ and $(s,i,-)$ cannot both occur while
calculating $T$. The reason for this is that the third component can
only be $+$ if either none of $x_1\cdots x_n$ has been read, in which
case the only relevant state is $(s_0,0,+)$, or else only $x_n$
has been read, in which case the possible relevant states are $(f,1,-)$,
$(s,1,+)$ with $s\neq f$, and $(s,0,-)$. So a state of the form
$(s,i,j)$ with a given $s$ occurs at most once in a fixed subset with the 
maximum possible value of $i$.

The effect of Step~\ref{modifieda} in the above algorithm is to ensure that
termination occurs as soon as a final state of $\Rev\Rules$ appears in a
calculated triple.
Since we know that $x_1\cdots x_n$ contains a left-hand
side of a rule in $\SetOfRules\Rules$ as a suffix we need only show that the 
introduction of Step~\ref{modifiedb} does not affect the accepted language of
the constructed automaton. This will be a consequence of
Theorem~\ref{modified determinization}, as we now proceed to show.

Consider a triple $t=(s,i,j)$ arising during the calculation of a subset
$T$, and suppose that $s$ is a non-final state of $\Rev\Rules$. 
If $j = +$ then $T$ cannot contain both $(s,0,+)$ and $(s,1,+)$ and so $t$ will
not be removed from $T$ as a result of Step~\ref{modifiedb}. Therefore we only
need to consider the case $j = -$. For $k=0,1,2$, let 
$L_k \subseteq \Astar\times\Astar$ be the language
obtained by making $(s,k,-)$ the only initial state of $M$, and observe that
there can be no padded arrows in any path of arrows from $(s,k,-)$ to a
final state of $M$. Now by considering the definition of the non-padded 
transitions in $M$ given in \ref{non-padded M transitions}, it is 
straightforward to see that $L_0 \subseteq L_1 = L_2$.
Therefore, since $\NRevLHS\Rules$ has no $\epsilon$-arrows, we have just 
shown that the hypotheses of Theorem~\ref{modified determinization} apply to
Step~\ref{modifiedb}. Hence the omission in
Step~\ref{modifiedb} does not affect the accepted language of $Q(\Rules)$. 
\end{proof}

As with $P(\Rules)$, reading a word into $Q(\Rules)$ from right to
left can be slow in the initial stages of a Knuth--Bendix pass, but soon
speeds up to being linear with a small constant.

\section{Finding the right-hand side of a rule}
\label{finding the right-hand side}

We retain the hypotheses of Section~\ref{properties of Rules}.
Namely, we have a two-variable rule automaton $\Rules$ which is welded
and satisfies various other minor conditions. We are given a word
$w=x_1\cdots x_n$, and we wish to reduce it relative to the rules implicitly
contained in $\Rules$.
So far we have located a left-hand side $\lambda$ which is
a substring of $w$. In this section we show how
to construct the corresponding right-hand side.

We first go into more detail as to how we propose to reduce $w$.
In outline we proceed as follows.

\subsection{Outline of the reduction process.}
\label{history stack}
\begin{enumerate}
\item Feed $w$ one symbol at a
time into the one-variable automaton $P(\Rules)$ described in
Section~\ref{which words are reducible},
storing the history of states reached on a stack.
\item If a final state is reached after some prefix $u$ of $w$ has been
read by $P(\Rules)$, then $u$ has some suffix which is a left-hand
side. Moreover, this procedure finds the shortest such prefix.
\item Feed $u$ from right to left into $Q(\Rules)$. A final state
is reached as soon as $Q(\Rules)$ has read the shortest suffix $\lambda$
of $u$ such that there is a rule $(\lambda,\rho) \in \SetOfRules\Rules$.
We now have $u = p\lambda$ and $w=p\lambda q$, where
$p,q \in \Astar$, every proper prefix of $p\lambda$
and every proper suffix of $\lambda$ is $\SetOfRules\Rules$-irreducible.
\item\label{use rule in Store}
Find $\rho$, the smallest string such that there
is a rule $(\lambda,\rho)$ in \Store (see \ref{Knuth--Bendix pass}).
If there is no such rule in \Store, find $\rho$ by a method to
be described in this section, such that $\rho$ is the smallest string
such that $(\lambda,\rho)\in \SetOfRules\Rules$.
\item\label{reduction gives new rule}
If $(\lambda,\rho)$ is not already in \Store, insert it into the
part of \Store called \New.
\item Replace $\lambda$ with $\rho$ in $w$ and pop $|\lambda|$ levels off
the stack so that the stack represents the history as it was immediately
after feeding
$p$ into $P(\Rules)$.
\item
Redefine $w$ to be $p\rho q$.
Restart at Step 1 as though $p$ has just been read and the next letter
to be read is the first letter of $\rho$. The history
stack enables one to do this.
\end{enumerate}

Note that other strategies might lead to
finding first some left-hand side in $w$ other than $\lambda$. Moreover,
there may be several different
right-hand sides $\rho$ with $(\lambda,\rho)\in
\SetOfRules\Rules$. A rule $(\lambda,\rho)$ in $\SetOfRules\Rules$ gives rise 
to paths in $\Rules$, $\SLtwo$ and $\DRev\SLtwo$. We will find the path for 
which right-hand side $\rho$ is \textit{shortlex}-least, given
that the left-hand side is equal to $\lambda$.

Let $\lambda = y_1\cdots y_{m}$.
Recall that a state of the one-variable automaton $Q(\Rules)$ used to find
$\lambda$ is a set of states of the form
$(s,i,j)$, where $s$ is a state of $\Rules, i\in\{0,1,2\}$ and $j\in\{+,-\}$.
When finding $\lambda$ we kept the history of states of $Q(\Rules)$ which were
visited---see Step~\ref{start}.
Let $Q_k$ be the set of triples $(s,i,j)$
comprising the state of $Q(\Rules)$ after reading
the string $ y_{k+1}\cdots y_m$ from right to left.
$Q_0=\{f_i\} = \{(s_0,i,-)\}$ where $s_0$
is the unique initial and final state of $\Rules$, and $i$ is the 
difference in length between $\lambda$ and the $\rho$ that we are looking for.

\subsection{Right-hand side routine.}
\label{right-hand side routine}
Inductively, after reading $y_1\cdots y_k$ we will have determined
$z_1\cdots z_k$, the prefix of $\rho$. Inductively we also have a triple
$(s_k,i_k,j_k)$, where $s$ is a state of $\Rules$, $i_k$ is 0 or 1 or
2 and $j_k$ is $+$ or $-$. Note that we always have $m-k \ge  i_k$.
\begin{enumerate}
\item\label{start2} If $m-k = i_k$, then we have found
$\rho = z_1\cdots z_k$ and we stop. So from now on we assume that
$m > i_k + k$. This means that the next symbol $(y_{k+1},z_{k+1})$ of
$(\lambda,\rho)$ does not have a padding symbol in its right-hand component.
\item We now try to find $z_{k+1}$ by running through each element $z\in A$ in
increasing order. Set $z$ equal to the least element of $A$.
\item \label{possible equality}
If $k=0$ and $i_0=0$, then $\lambda$ and $\rho$ will be of
equal length, so the first symbol of $(\lambda,\rho)$ must be
$(y_1,z_1)$, where $y_1> z_1$.
So at this stage we can prove that we have $y_1 > z$, since we know
that there must be some right-hand side corresponding to our given
left-hand side.

If $k=0$ and $i_0 > 0$, then the first symbol of
$(\lambda,\rho)^+$ is $(y_1,z_1)$ with $z_1\in A$ and $y_1\neq z_1$.
If $k=0$, $i_0>0$ and $y_1=z$, we increase $z$ to the next element of $A$.
\item\label{find z}
Here we are trying out a particular value of $z$
to see whether it allows us to get further.
We look in $\Rules$ to see if $s_k^{(y_{k+1},z)} = s_{k+1}$ is defined.
If it is not defined, we increase $z$ to the next element of $A$
and go to Step~\ref{possible equality}.

\item\label{finding the triple}
 If $s_{k+1}$ is defined in Step~\ref{find z}, we look in $Q_{k+1}$ for a
triple $(s_{k+1},i_{k+1},j_{k+1})$ which is the source of an arrow
labelled $(y_{k+1},z)$ in the automaton $M$.
Recall that $M$ was defined in
Section~\ref{finding the left-hand side}.
Note that, by the proof of \thmref{lhs algorithm}, $Q_{k+1}$ contains at most 
one element whose first coordinate is $s_{k+1}$. As a result, the search can be
quick.
\item If $(s_{k+1},i_{k+1},j_{k+1})$ is not found in Step~\ref{finding the triple},
increase $z$ to the next element of $A$ and go to Step~\ref{possible equality}.
\item If $(s_{k+1},i_{k+1},j_{k+1})$ is found in Step~\ref{finding the triple},
set $z_{k+1} = z$, increase $k$ and go to Step~\ref{start2}.
\end{enumerate}

The above algorithm will not hang, because each triple $(s_k,i_k,j_k)$
that we use does come from a path of arrows in $M$ which starts
at the initial state of $M$ and ends at the first possible final state
of $M$. Therefore all possible right-hand sides $\rho$ such that
$(\lambda,\rho) \in \SetOfRules\Rules$, are implicitly computed when we
record the states of $Q(\Rules)$ (see
Step~\ref{start}).
Since $i_k$ does not vary during our search, we
will always find the shortest possible $\rho$, with $|\lambda| -
|\rho| $ being equal to this constant value of $i_k$.
Since we always look for $z$ in increasing order, we are
bound to find the lexicographically least $\rho$.

We remind the reader that an overview of the entire reduction process
for a given string $w$ is given in \ref{history stack}. 

\section{Our version of Knuth--Bendix.}
\label{our version}
For finite Noetherian rewriting systems the question of confluence is
decidable by the critical pair analysis described in Section~\ref{strrewriting}.
However, for
{\it infinite} Noetherian rewriting systems the confluence question is, in
general,
undecidable. Examples exhibiting undecidability are given in 
\cite{undecidability} and are length-reducing rewriting
systems $R$ which are {\it regular} in the sense that $R$ contains only a finite
number of right-hand sides and for each right-hand side $r$, the set
$\{l : (l,r) \in R\}$ is a regular language.
These examples are in the context of rewriting for monoids.
As far as we know,
there is no known example of undecidability if we add to the hypothesis
that the monoid defined by $R$ is in fact a group.

In this section we consider a
rewriting system which is the accepted language of a rule automaton for some
finitely presented group. We  describe a Knuth--Bendix type algorithm for such
a system. In light of the undecidability result mentioned above, our algorithm
does not provide a test for confluence.
We can however use our algorithm together with other algorithms for
dealing with \textit{shortlex}-automatic groups,
to prove confluence by an indirect
route if the group is \textit{shortlex}-automatic.
Details of the theory of how this is done can be found in
\cite{wordprocessing}. The practical details are carried out in
programs by Derek Holt---see \cite{KBMAG}.

Suppose throughout that $G$ is a monoid given by a finite
presentation $\langle A / R \rangle$, where $A$ is a
set of generators for $G$ with a fixed total ordering
$<$ and $R$ is a finite set of equalities. The monoid is defined by the
congruence generated by these equalities.
We will assume that there is an involution $\iota:A\to A$ (which will
send each generator to its formal inverse) such that,
for each element $x\in A$, there are equalities in $R$ of the form
$x\iota(x)=\epsilon$ and $\iota(x)x=\epsilon$.
This implies that $G$ is a group.
The equalities in $R$ can be regarded
as a finite set of rules which define $G$.

In our algorithm, we keep two sets of rules. One set, which we call \Store,
is a finite set of rules. The other is a possibly infinite set of rules
which is kept implicitly in a rule automaton called $\Rules$. When we want to
specify that we are working with the $\Rules$ automaton during the $n$th 
Knuth-Bendix pass (see \ref{Knuth--Bendix pass} for the definition of a
Knuth-Bendix pass), we will use the notation $\Rules[n]$. We extract
explicit rules from $\Rules[n]$ by taking elements of the intersection
$\SetOfRules{\Rules[n]} = L(\Rules[n]) \cap L(\SLtwo)$. The two-variable
automaton $\SLtwo$ was defined in Section~\ref{which words are reducible} and
is depicted in \figref{SL2}. 

\Store will change almost continually, while $Rules$ is constant during a
Knuth-Bendix pass. We do in fact need to change $\Rules$ from time to
time, and we do so as the last step of each Knuth--Bendix pass.
We will perform the Knuth--Bendix process, using the rules in \Store
for critical pair analysis, as described in \ref{critical case analysis}.

\subsection{Rapid reduction.}
\label{rapid reduction}
A difference between our situation and that of
classical Knuth--Bendix is that reduction is not
carried out by applying the rules of \Store.
When running Knuth--Bendix, one of the most time-consuming aspects is
reduction. This is partly because there is a lot of reduction to be
done and partly because one normally has to spend a long time
looking through a long list of rules to see if the string one is
trying to reduce contains a left-hand side of some rule.
Much of the effort in producing new Knuth--Bendix algorithms, like the
algorithm described in this paper,
goes on finding methods of locating relevant rules quickly.
In the past this has involved using structures which use a lot of
space. In our procedure we use the method described in \ref{history stack} to 
find relevant rules quickly without using an inordinate amount of space.
We refer to this as \textit{\RR-reduction}. We also use the terms
\textit{\RR-reduce} and its various derivatives.
\RR stands for ``relation'', for ``reduction'' and for ``rapid''.

\begin{note}\label{changing RR}
Note that a string is \RR-reducible at one point in a Knuth--Bendix
pass if and only if it is \RR-reducible at another point in the same
pass. However, as we shall see,
the result of \RR-reduction may change during a pass,
because \Store changes.
\end{note}

\subsection{The basic structures.}
The basic structures used in our procedure are:
\begin{enumerate}
\item A two-variable automaton $\Rules$ satisfying the conditions
laid down in \ref{properties of Rules}.
\item A finite set \Store of rules, which is the disjoint union of
several subsets of rules : \Considered, \This, \New and \Delete.
\item \Considered is a subset of \Store such that each rule
has already been compared with each other rule in \Considered,
including with itself,
to see whether left-hand sides overlap.
The consequent critical pair
analysis has also been carried out for pairs of rules in \Considered.
Such rules do not need to be compared with
each other again.
\item \This is a subset of \Store containing rules which we plan to use
during this pass to compare for overlaps with the rules in \Considered, as in
\ref{first case}. These rules have been
minimized during the current pass (see \ref{minimal}) and so should not
be minimized again.
\item \New is a subset of \Store containing new rules which have
been found during the current pass, other than those which are output
by the minimization routine (see \ref{minimal}). Non-trivial rules which are the
final output of the minimization routine are added to \This.
\item \Delete is a subset of \Store containing rules which are to be
deleted at the beginning of the next pass.
\item A two-variable automaton \WDiff which contains
all the states and arrows of $\Rules[n]$, and possibly other states and
arrows. It satisfies the conditions of \ref{properties of Rules}.
\end{enumerate}

\subsection{Initial arrangements.}
Before describing the main Knuth--Bendix process, we explain how the
data structures are initially set up. Recall that $R$ (which should be distinguished
from \RR) is the original set of
defining relations together with special
rules of the form $(x\iota(x),\epsilon)$ and $(\iota(x)x,\epsilon)$ which
make the formal inverse $\iota(x)$ into the actual inverse of $x$.

We rewrite each non-special element of $R$ in the form
of a relator, which we cyclically reduce in the free group.
Since $\iota(x)$ is the formal inverse of the letter $x$, we are able to
write down the formal inverse of any string in $\Astar$.
We may therefore assume that each relator has the form
$lr^{-1}$, where $l$ and $r$ are elements of $\Astar$ and $(l,r)$ is
accepted by $\SLtwo$.

For each rule $(l,r)$, including the special rules, we form a rule
automaton, as explained in \exref{ruleworddiff}. These automata are then welded
together to form the two-variable rule automaton \WDiff satisfying the 
conditions of \ref{properties of Rules}. Each state and arrow 
of \WDiff is marked as \needed. (At certain well-chosen moments we will
delete from \WDiff states and arrows that are not \needed). Each of these rules
is inserted into \This. \Considered, \New and 
\Delete are initially empty. Set $\Rules[1]=\WDiff$.

\subsection{The main loop---a Knuth--Bendix pass.}
A significant proportion of the time in a Knuth--Bendix pass is spent in
applying a procedure which we term {\it minimization}. Each rule encountered
during the pass is input to this procedure and the output is called a 
{\it minimal rule}. 
The exact details of this process are given in sections~\ref{minimizing a rule}
and \ref{handling minimization output}, but we point out that minimization 
often results in rules being added to and/or deleted from \Store. Any rules
added to \Store during the minimization of a rule $(\lambda,\rho)$ are
strictly smaller than $(\lambda,\rho)$ in the ordering of 
\ref{ordering details}. 

\label{main loop}
\begin{enumerate}
\item\label{start main loop}
At this point, \This is empty.
If $n>0$, save space by deleting previously defined automata $P(\Rules[n])$,
$Q(\Rules[n])$ and $\Rules[n]$. Increment $n$.
The integer $n$ records which Knuth--Bendix pass we are currently working on.
\item \label{process Delete}
Delete the rules in \Delete.
\item\label{process Considered}
For each rule $(\lambda,\rho)$ in \Considered,
minimize $(\lambda,\rho)$ as in \ref{minimal} and handle the output as in
	\ref{handling minimization output}.
\item \label{process New}
For each rule $(\lambda,\rho)$ in \New,
minimize $(\lambda,\rho)$ as in \ref{minimal} and handle the output as in
	\ref{handling minimization output}.

Since rules added to \New during minimization are always strictly smaller than
than the rule being minimized, it follows that eventually each rule in \New 
will be processed; that is, the process of examining rules in \New does not 
continue indefinitely.
\item\label{process This}
For each rule $(\lambda,\rho)$ in \This:
	\begin{enumerate}
	\item Delete the rule from \This and add it to \Considered.
	\item\label{compare for overlaps}
		 For each rule $(\lambda_1,\rho_1)$ in \Considered:
		\begin{quote} Look for overlaps between $\lambda$ and
		$\lambda_1$. That is we have to find each
		suffix of $\lambda$ which is a prefix
		of $\lambda_1$ and
		each suffix of $\lambda_1$ which is a
		prefix of $\lambda$.
		Note that we may have to allow $\lambda=\lambda_1$ in order to
		deal with the case where two different rules have the
		same left-hand side. In this case, both the prefix and suffix 
		of both left-hand sides is equal to $\lambda=\lambda_1$.
		Then \RR-reduce in two different ways as in
		\ref{first case}, obtaining a pair of strings $(u,v)$.
		If they differ then rearrange them so that $u>v$ and
		insert the result into \New, unless it is already in
		\Store.
		\end{quote}
	\end{enumerate}
\item\label{using WDiff}
Delete from \WDiff all arrows and states which are not marked as
\needed. Copy \WDiff into $\Rules[n+1]$
and mark all arrows and states of \WDiff as not \needed.
The details of this step are given in \ref{WDiffdetails}.
\item This ends the description of a Knuth--Bendix pass.
Now we decide whether to terminate the Knuth--Bendix process. Since we know of 
no procedure to decide confluence of an infinite system of rules (indeed, it is
probably undecidable), this decision is taken on heuristic grounds. In our
context, a decision to terminate could be taken simply on the grounds that 
\WDiff and $\Rules[n]$ have the same states and arrows. In other words, no new
word-differences or arrows between word-differences has been found
during this pass.
If the Knuth--Bendix process is not terminated, go to \ref{start main loop}.
\end{enumerate}

\subsection{Minimizing a rule.}
\label{minimizing a rule}
We now provide the details of the minimization routine. 
\begin{defn}
\label{minimal}
Let $(u,v)\in \Astar\times\Astar$ and let $u=u_1\cdots u_p$
and $v=v_1\cdots v_q$, where $u_i, v_j \in A$.
We say that $(u,v)$ is a \textit{minimal rule} if $u\neq v$,
$\bar u = \bar v$ in $G$
and the following procedure does not change $(u,v)$. The procedure is
called \textit{minimizing a rule} or the \textit{minimization routine}.
We always start the minimization routine with $u>_{SL} v$, though this
condition is not necessarily maintained as $u$ and $v$ change during the
routine.
\begin{enumerate}
\item \RR-reduce the maximal proper prefix $u_1\cdots u_{p-1}$
of $u$ obtaining $u'$.
Reduction may result in rules being added to \New as described in
\ref{reduction gives new rule}.
If $u\neq u'u_p$, change $u$ to $u'u_p$ and go to
Step~\ref{minimize u}.
\item \RR-reduce the maximal proper suffix $u_2\cdots u_p$
of $u$ obtaining $u''$.
Reduction may result in new rules being added to \New.
Replace $u$ by $u_1u''$.
\item\label{minimize u}
If $u$ has changed since the original input to the minimization routine,
then \RR-reduce $u$. This may result in new rules being added to \New.
\item\label{minimization loop}
If $p >q+2$ or if $p = q+2$ and $u_1 > v_1$,  replace $(u,v)$ by
$(u_1\cdots u_{p-1}, v_1\cdots v_q \iota(u_p))$ and repeat this step until
we can go no further.
\item If $p=q+2$ and $u_2 > \iota(u_1)$, replace $(u,v)$ by
$(u_2\cdots u_p,\iota(u_1)v_1\cdots v_q)$.
\item If $q>0$ and $u_1 = v_1$, cancel the first letter
from $u$ and from $v$ and repeat this step.
\item If $q>0$ and $u_p =v_q$, cancel the last letter from $u$ and
from $v$ and repeat this step.
\item \RR-reduce $v$ as explained in \ref{history stack}.
This may result in rules being added to \New as described in
\ref{reduction gives new rule}.
\item If $v > u$, interchange $u$ and $v$.
\item If $(u,v)$ has changed since the last time Step~\ref{minimization
loop} was executed, go to Step~\ref{minimization loop}.
\item Output $(u,v)$ and stop.
\end{enumerate}
From Note~\ref{changing RR}, we see that if a rule is minimal at
one time during a Knuth--Bendix pass then it is minimal at
all later times during the same pass.
\end{defn}

Note that the output could be $(\epsilon,\epsilon)$, which means
that the rule is
redundant. Otherwise we have output $(u,v)$ with $u>v$.
Note that the minimization procedure keeps on decreasing $(u,v)$ in the
ordering given by \ref{ordering details}.
Since this is a well-ordering, the minimization
procedure has to stop.
Also any rules added to \New as a result of \RR-reduction during
minimization are smaller than $(u,v)$.

\begin{lemma}\label{proper substrings irreducible}
Let $(\lambda_1,\rho_1)$ be the output from minimizing $(\lambda,\rho)$.
If $\lambda$ has no proper \RR-reducible substrings, then
$\lambda_1$ is a non-trivial substring of $\lambda$.
\end{lemma}
\begin{proof}
Under the hypotheses, the successive steps of
minimization change $\lambda$ and $\rho$, while maintaining
the inequality $\lambda > \rho$.
As a result the left-hand and right-hand sides are never interchanged.
It follows that $\lambda_1 > \rho_1$, so $\lambda_1$ is
non-trivial. It is easy to see that $\lambda_1$ is a substring of $\lambda$.
\end{proof}

\subsection{Handling minimization output.}
\label{handling minimization output}
Suppose the input to minimization is $(\lambda,\rho)$ and its output
is $(\lambda_1,\rho_1)$. We now describe how the rule $(\lambda_1,\rho_1)$ is
treated. 
\begin{enumerate}
\item In order to avoid unnecessary subsequent work of minimization,
mark $(\lambda,\rho)$ as minimized and $(\lambda_1,\rho_1)$ as
minimal for this pass.
\item If $(\lambda_1,\rho_1) \neq (\epsilon,\epsilon)$,
incorporate $(\lambda_1,\rho_1)$ into the language accepted
by \WDiff, using the method described in \ref{WDiffdetails}.
\item If $(\lambda_1,\rho_1)=(\lambda,\rho)$, that is, if $(\lambda,\rho)$
was already minimal, then:
	\begin{enumerate} \item If $(\lambda,\rho)$ was in \Considered, do
	nothing.
	\item If $(\lambda,\rho)$ was in \New, move it to \This.
	\item We saw in \ref{main loop}, that the above two situations
	are the only ones in which minimization is carried out.
	\end{enumerate}
\item If $(\lambda,\rho)\neq(\lambda_1,\rho_1)\neq (\epsilon,\epsilon)$, then:
	\begin{enumerate}\item If $(\lambda_1,\rho_1)$ is already in
	\Store, mark the copy in \Store as minimal.
	\item If $(\lambda_1,\rho_1)$ is not already in \Store, insert it
	in \This.
	\end{enumerate}
\item If $(\lambda_1,\rho_1) = (\epsilon,\epsilon)$, do nothing.
\item\label{deletion: proper substring} If $\lambda$ was affected by the minimization routine before
Step~\ref{minimization loop}, that is, if some proper substring of $\lambda$
was \RR-reducible, then delete $(\lambda,\rho)$.
\item If, at the time of minimization,
all proper substrings of $\lambda$ were \RR-irreducible and
if $(\lambda,\rho)$ was not minimal, move
$(\lambda,\rho)$ to the \Delete list. The reason for this possibly
surprising policy of not deleting immediately is that further
reduction during this pass may once again
produce $\lambda$ as a left-hand side by the
methods of Sections~\ref{which words are reducible} and
\ref{finding the left-hand side}. We want to avoid the work involved in
finding the right-hand side by the method of
Section~\ref{finding the right-hand side}. For this, we need to have a rule
in \Store with left-hand side equal to $\lambda$---see \ref{reduction
gives new rule}.
\end{enumerate}

\subsection{Details on the structure of \WDiff.}
\label{WDiffdetails}
At the beginning of Step~\ref{using WDiff}, each state $s$ of $\WDiff$ has an
associated string $w_s \in A\uast$ which is irreducible with respect
to $\SetOfRules{\Rules[n]}$. $\WDiff$ is a rule automaton: we associate the
element $\overline{w_s}$ to the state $s$.
These state labels are calculated as and when new states and arrows are added 
to $\WDiff$ during a Knuth--Bendix pass (see \ref{sewing}).

At the end of the $n$th Knuth--Bendix pass,
\WDiff is an automaton which represents the
word-differences and arrows between them encountered during that pass.
At this stage the string
attached to each state is irreducible with respect to the rules in
$\SetOfRules{\Rules[n]}$ but not necessarily 
with respect to the rules implicitly contained in \WDiff. Before starting the
next pass, we \RR-reduce the state labels of \WDiff with respect to 
$\SetOfRules{\WDiff}$. If \WDiff now contains distinct states labelled by the
same string we connect them by epsilon arrows and replace \WDiff by
$\Weld \WDiff$ (see remark~\ref{identify}). We then repeat this procedure until
all states are labelled
by distinct strings which are irreducible with respect to 
$\SetOfRules{\WDiff}$. If during this procedure a state or arrow marked
as \needed is identified with one not marked as \needed, the resulting
state or arrow is marked as \needed.

Whenever a minimal rule $r$ is encountered during the $n$th pass, it 
is adjoined to the accepted language of \WDiff. One method of doing this is
to form the rule automaton $M(r)$ as given in example~\ref{ruleworddiff}, and
replace \WDiff by the result of welding the union of the two rule automata
\WDiff and $M(r)$. However, there is a more efficient way to proceed which
eliminates the need to construct $M(r)$ and possibly avoids the necessity
for welding. We call this procedure
{\bf sewing}. 

\subsection{The sewing procedure.}
\label{sewing}
Suppose we have a rule $r=(u_1,v_1)\cdots (u_n,v_n)$, with each
$u_i\in A$ and each $v_i\in A^+$ which we want to add to the language
accepted by \WDiff.
We read the rule into \WDiff from left to right, starting at the initial state of
\WDiff. Suppose it is possible to read in $(u_1,v_1) \cdots (u_k,v_k)$,
arriving at a state $s_k$, where $k$ is chosen as large as possible,
subject to the condition that $k\le n$.
If $k<n$, then the arrow labelled $(u_{k+1},v_{k+1})$
is undefined on $s_k$.
We now read $(u_{k+1},v_{k+1}) \cdots (u_n,v_n)$ into \WDiff, starting from
the initial and final state $s_0 = t_n$ and reading from
right to left, arriving at a state $t_r$ with $r\ge k$,
on which the backward arrow
labelled $(u_r,v_r)$ is undefined.
We mark all states and arrows encountered as \needed.

We now proceed as follows:
\begin{enumerate}
\item \label{creating}
If $k=r$ and the states $s_k$ and $t_r$ do not coincide, join them by an
epsilon-arrow, replace \WDiff by $\Weld\WDiff$, and then stop. 
When identifying states with different labels during the welding procedure, we
choose the {\it shortlex}-least label for the amalgamated state.
\item If $k<r$, let $w_k$ be the label for $s_k$. Reduce $u_{k+1}^{-1}w_k
v_{k+1}$, obtaining $w_{k+1}$.
\item If $w_{k+1}$ is the label of an existing
state $t$ of \WDiff, set $s_{k+1} = t$.
\item If $w_{k+1}$ is not the label of an existing state, create a new state
$s_{k+1}$ with label $w_{k+1}$.
\item Create an arrow labelled $(u_{k+1},v_{k+1})$ from $s_k$ to $s_{k+1}$.
\item Mark the state $w_{k+1}$ and the arrow $(u_{k+1},v_{k+1})$ as
\needed.
\item Increment $k$ and go to Step~\ref{creating}.
\end{enumerate} 

Note that the automaton obtained by sewing is a rule automaton.

\section{Correctness of our Knuth--Bendix Procedure}
\label{correctness}
In this section we will prove that the procedure set out in
Section~\ref{our version} does what we expect it to do.

\begin{defn} For a discrete time $t$, we denote by \Storen t the rules in
\Store at time $t$ in our Knuth---Bendix procedure. We can take $t$ to be
the number of machine operations since the program started, or any similar
discrete measure.
\end{defn}

\begin{defn}
A quintuple $(t,s_1,s_2,\lambda,\rho)$, where $t$ is a time, and $s_1$,
$s_2$, $\lambda$ and $\rho$ are elements of $\Astar$, is called an
{\it elementary \Storen t-reduction} $u \rightarrow_{\Storen t} v$ from $u$
to $v$ if $(\lambda,\rho)$ is a rule in \Storen t, $u=s_1\lambda s_2$
and $v=s_1\rho s_2$. We call $(\lambda,\rho)$ the \textit{rule associated
to the elementary reduction}.
\end{defn}

\begin{defn}\label{Thue path}
Let $t\ge 0$.
By a \textit{time-$t$ Thue path} between two strings $w_1$ and $w_2$,
we mean a sequence of elementary \Storen t-reductions and inverses of elementary
\Storen t-reductions connecting $w_1$ to $w_2$, such that none of the
rules associated to the elementary reductions is in \Delete at time $t$.
We talk of the strings which are the source or target of these elementary
reductions as \textit{nodes}.
The path is considered as having a direction from $w_1$ to $w_2$.
The elementary reductions will be consistent with this
direction and will be called \textit{rightward} elementary reductions.
The inverses of elementary reductions will be in the opposite direction
and will be called \textit{leftward} elementary reductions.
\end{defn}

\begin{proposition}
\label{consistent rules}
Let $\langle A / R \rangle$ be the finite presentation at the start of
the Knuth-Bendix procedure. Then if $(\lambda,\rho)\in\Store$ during the $n$-th
Knuth-Bendix pass, we have $\lambda \leftrightarrow_R^* \rho$.
\end{proposition}

\begin{proof}
The proof of this is an easy induction on $n$ using 
Corollary~\ref{validrulescorollary}.
\end{proof}

\begin{proposition}
\label{maintain congruence}
Let $t\geq0$ and suppose that we have a Thue path from $\alpha$ to 
$\beta$ in \Storen t with maximum node $w$. Then for any time $s\geq t$,
there exists a time-$s$ Thue path from $\alpha$ to $\beta$
with each node less than or equal to $w$.
\end{proposition}
\begin{proof}
We show by induction on $s$ that, if at some time $t\leq s$ there is a Thue
path between strings $u$ and $v$ with all nodes no bigger than $u$ or $v$, then
there is also such a Thue path at time $s$. So suppose that we have proved this
statement for all times $s'<s$. 

We first consider the special case where $(u,v)$ is
a rule being input to the minimization routine (see Definition~\ref{minimal}) 
at time $t$, and $s$ is the time at end of the subsequent invocation of the
minimization handling routine \ref{handling minimization output}. 

Each step of minimization takes an input string and outputs
a possibly different string which is used as the input to the next step. The
initial input is $(u,v)$ and the final output is either $(\epsilon,\epsilon)$
or a minimal rule $(u',v')$. Let $r_1,r_2,\ldots,r_n$ be the sequence
of outputs in the minimization of $(u,v)$, and let $r_0=(u,v)$. By
considering each step of minimization in turn, we will show that for each
$i,1\leq i\leq n$, if there is a time-$s$ Thue path between the two sides of
$r_i$ with maximum node no bigger than either side, then there is a time-$s$
Thue path between the two sides of $r_{i-1}$ with maximum node no bigger than
either side. We then obtain the desired time-$s$ Thue path between $u$ and
$v$ by using descending induction on $i$ given that the base case $i=n$ is 
trivially true. 

To make the task of checking the proof easier, we use the same numbering
here as in Definition~\ref{minimal}.
\begin{enumerate}
\item At the end of this step, there is a sequence of elementary reductions
from $u_1\ldots u_{p-1}$ to $u'$, but this may not constitute a Thue path
since some of the associated rules may be in \Delete. However, any such rule 
$(\lambda,\rho)$ will, at some time $s' < s$, have been in \Store but not in 
\Delete. Therefore by our induction on $s$, at the end of this step there will 
be a Thue path from $\lambda$ to $\rho$ with maximum node $\lambda$. Since
no rule used in this Thue path is equal to $(u,v)$, this will still be a
Thue path at time $s$. Hence we can construct a time-$s$ Thue path from $u$ to
$v$ with maximum node $u$. 

\item This step is analogous to the previous step.

\item At the end of this step, the sequence of \RR-reductions of $u$ to the 
current left-hand side does not use the rule $(u,v)$ (hence the condition at
the start of this step), and so the required Thue path exists. 

\item Suppose that the input to this step is $(u'x,v)$.
Then the output is either the same as the input or is equal to
$(u',vx^{-1})$. In the first case there is nothing to prove. In the latter 
case, a time-$s$ Thue path from $u'$ to 
$vx^{-1}$ with maximum node $u'$ will give a time-$s$ Thue path from $u'x$
to $vx^{-1}x$ with maximum node $u'x$. Note that there have been no deletions of
rules since this particular
minimization was started. Induction on $s$ therefore gives us a specific Thue 
path from $x^{-1}x$ to $\epsilon$ at the end of this step.
Moreover, no node along the Thue path is
bigger than $x^{-1}x$. In particular, the input rule to the minimization
is not used in this Thue path (this follows from either of the conditions at 
the start of the step 
being satisfied), and so there is a time-$s$ Thue path from $vx^{-1}x$ to $v$ 
with maximum node $vx^{-1}x$. Hence we obtain the required time-$s$ Thue path
from $u'x$ to $v$. 
\item This step is analogous to the previous step. 
\item If the input to this step is $(xu',xv')$ then the output is $(u',v')$.
A time-$s$ Thue path from $u'$ to $v'$ with maximum node $u'$ yields
a time-$s$ Thue path from $xu'$ to $xv'$ with maximum node $xu'$.
\item This step is analogous to the previous step. 
\item Let $v'$ be the \RR-reduction of $v$. Immediately after this step there
is a Thue path from $v$ to $v'$ with maximum node $v$ which does not use the 
rule initially input.
Using induction if necessary we have a
time-$s$ Thue path from $v$ to $v'$ with maximum node $v$. Hence a time-$s$
Thue path from $u$ to $v'$ with maximum node either $u$ or $v'$ yields a 
time-$s$ Thue path from $u$ to $v$ with maximum node either $u$ or $v$.
\item If there is a Thue path from $u$ to $v$ with maximum node either
$u$ or $v$, then the reverse of this path is a Thue path from $v$ to $u$.
\end{enumerate}
This completes the induction step for the special case. Now consider the
general case. The only reason why a Thue path at time $t < s$ between $u$ and
$v$ will not work at time $s$ is if some elementary reduction used in this
path has an associated rule $(\lambda,\rho)$ in \Storen t which is found to be
non-minimal between $t$ and $s$. But in the proof of the special case we have
seen that there is a time-$s$ Thue path between $\lambda$ and $\rho$ with
no node bigger than $\lambda$. Therefore the time-$t$ Thue path can always be
replaced by a time-$s$ Thue path without increasing the size of the nodes. 
\end{proof}

\begin{lemma}\label{stays S-reducible}
If a string is \Storen s-reducible, it is \Storen t-reducible for all
$t>s$.
\end{lemma}
\begin{proof}
If $u$ is \Storen s-reducible, there is an elementary
\Storen s-reduction $u\to_{\Storen s} v$.
This means that $v<u$. By Proposition~\ref{maintain congruence},
for each time $t>s$, there is a Thue path from $u$ to $v$ with maximum
node $u$. The first elementary reduction in this path has the form
$u\to w$ at time $t$. This proves the result.
\end{proof}

\begin{lemma}\label{stays R-reducible}
If $(\lambda,\rho)$ is a rule in \Store at some time during
the $n$-th Knuth--Bendix pass but before the beginning of 
Step~\ref{process This}, then $\lambda$ will be \RR-reducible during all 
subsequent passes. If $\lambda$ is \RRn s-reducible then $\lambda$ is \RRn 
t-reducible for any $t>s$.
\end{lemma}
\begin{proof}
Let $(\lambda,\rho)$ be a rule as in the statement of the lemma. Then at some
prior time, $(\lambda,\rho)$ will have been a rule in \Store but not in
\Delete. Therefore for any $m\geq n$, there will be a Thue path from 
$\lambda$ to $\rho$ with maximum node $\lambda$ at the beginning of 
Step~\ref{process This} during the $m$-th Knuth-Bendix pass. Now at the 
beginning of Step~\ref{process This}, all rules in \Store but not in \Delete 
will have been output by the minimization handling routine 
\ref{handling minimization output} at some prior time during that Knuth-Bendix
pass. In particular, each of these minimal rules $(u,v)$ will have been sewn 
into \WDiff. This does not, however, imply that $(u,v)$ will be accepted by
$\Rules$ at the start of the next pass since this rule may use or define an
$(x,x)$ arrow in \WDiff. Due to some collapsing in \WDiff caused by a
welding operation, this may give rise to an $(x,x)$ arrow from the initial
state to itself. Such an arrow will be removed so that \WDiff satisfies the
properties \ref{properties of Rules}. If this is the case then $(u,v)$ will
still have some prefix or suffix accepted by \WDiff and hence by \Rules at
the start of the next Knuth-Bendix pass. Therefore for any $m\geq n$,
$\lambda$ will have a substring which is the left-hand side of a rule
accepted by $\Rules[m+1]$, and so $\lambda$ is \RR-reducible during pass
$m+1$ which proves the first statement.

If $\lambda$ is \RR-reducible at any time during a pass, it is \RR-reducible
at any later time in the same pass by Lemma~\ref{changing RR}.
We have proved that $\lambda$ is \RR-reducible once the next pass starts.
So this completes the proof of the last sentence in the statement of the
lemma.
\end{proof}

\begin{lemma}\label{no repeats}
At any time $t$, $\Storen t$ contains no duplicates.
If a rule is deleted from \Store, it 
will never be re-inserted.
\end{lemma}
\begin{proof}
The first statement follows by looking through \ref{main loop} and
checking where insertions of rules take place. We always take care not
to insert a rule a second time if it is already present.

Let $(\lambda,\rho)$ be a rule which is deleted at time $s$.
Deletion either takes place during Step~\ref{process Delete}
or during Step~\ref{deletion: proper substring}.
In the latter case,
some proper substring of $\lambda$ is the left-hand side of a rule in \Store,
and this rule was present in \Store at some time before the beginning of 
Step~\ref{process This} of the Knuth-Bendix pass in which $(\lambda,\rho)$ was
deleted. Therefore by Lemma~\ref{stays R-reducible}, we see that this
proper substring of $\lambda$
stays \RR-reducible. This means that no rule with left-hand side $\lambda$
will ever be re-inserted.

So we assume that deletion takes place during Step~\ref{process Delete} of
the $(n+1)$-th Knuth-Bendix pass. Then between the time of the $n$-th Knuth-
next pass when $(\lambda,\rho)$ is minimized, and the start of the next pass,
it is in the subset \Delete of \Store. Therefore it cannot be re-inserted 
during pass $n$. 

Now let $m>n$ and suppose $(\lambda,\rho)$ has not been re-inserted before
the beginning of the $m$-th pass.
We will prove that it cannot be re-inserted during the $m$-th pass.

Observe that no rules are minimized between the time $r$ at the beginning of
Step~\ref{process This} in the $(m-1)$-st pass and the time $t$ just
defined. Therefore any time-$r$ Thue path between $\lambda$ and $\rho$ will
also be a time-$t$ Thue path. In particular, the rule 
$(\lambda',\rho')$ associated to the elementary reduction of $\lambda$ is
unaltered during this time.
At time-$r$ all rules in \Store, except for those in \Delete, are minimal.
In particular, $(\lambda',\rho')$ was minimal
and had been minimized at some prior point in the $(m-1)$-st pass. Therefore,
$(\lambda',\rho')$ was sewn into \WDiff during the $(m-1)$-st pass. As in
the proof of Lemma~\ref{stays R-reducible}, it is possible that
$(\lambda',\rho')$ is not
accepted by $\Rules[m]$, but a substring $\lambda^{\prime\prime}$ of 
$\lambda'$ (and hence $\lambda$), will be the left-hand side of an accepted 
rule. If $\lambda^{\prime\prime}$ is a proper substring of $\lambda$, then 
$\lambda$ cannot be the left-hand side of a rule inserted during the $m$-th
Knuth--Bendix pass. So we need only examine the case when
$\lambda=\lambda'=\lambda^{\prime\prime}$. In this situation, $\rho$ must have
been \RR-reduced to a strictly smaller string during the minimization of 
$\lambda$.
Therefore some substring of $\rho$ was the left-hand side of a rule in \Store 
at that time. By Lemma~\ref{stays R-reducible}, $\rho$ stays \RR-reducible.

Suppose then that $(\lambda,\rho)$ is re-inserted during the $m$-th 
Knuth--Bendix pass. This can only happen as the result of
Step~\ref{reduction gives new rule}. But for this to occur there could have been
no rule in \Store with left-hand side equal to $\lambda$ at that time. 
Since we are assuming that $\lambda$ is the left-hand side of some rule
$(\lambda,\rho^{\prime\prime})$ not on the \Delete list at the start of the
$m$-th pass, it follows that $(\lambda,\rho^{\prime\prime})$ must be deleted at
some point during the $m$-th pass. But this can only happen if some proper 
substring of $\lambda$ is \RR-reducible during the $m$-th pass. By 
Lemma~\ref{changing RR}, this proper substring of $\lambda$ must be 
\RR-reducible at the point of re-insertion of $(\lambda,\rho)$ which is a
contradiction. 
\end{proof}

\begin{defn} We say that a string $u$ is \textit{permanently
irreducible} if there are arbitrarily large times $t$ for which
$u$ is \Storen t-irreducible. By Lemma~\ref{stays S-reducible}
this is equivalent to saying
that $u$ is \Storen t-irreducible at all times $t\ge 0$.
A rule $(\lambda,\rho)$ in \Store
is said to be \textit{permanent} if $\rho$ and
every proper substring of $\lambda$ is permanently irreducible.
\end{defn}

\begin{lemma}\label{never deleted}
A permanently irreducible string is permanently \RR-irreducible.
A permanent rule of \Store is never deleted.
A permanent rule is accepted by $\Rules[n+1]$ provided it is present in
\Store when the $n$-th Knuth--Bendix pass begins, {\rm(}and is accepted
by $\Rules[m]$ for all $m > n${\rm)}.
\end{lemma}
\begin{proof}
Let $u$ be permanently irreducible.
\RR-reduction of $u$ can only take place if, immediately afterwards, some
substring of $u$ is \Store-reducible. This is impossible by
hypothesis.

A rule $(\lambda,\rho)$ is deleted only as a result of minimization.
By Lemma~\ref{maintain congruence}, there would have to be a Thue path
from $\lambda$ to $\rho$ with largest node $\lambda$. The first
elementary reduction must therefore be rightward (see
Definition~\ref{Thue path}) $\lambda\rightarrow_{\Storen t} \mu$. Since every 
proper substring of $\lambda$ is permanently \RR-irreducible, this first
elementary reduction must be associated to a rule $(\lambda,\mu)$.

This is only possible if, when $(\lambda,\rho)$ was input to the minimization
routine, $\rho$ was \RR-reducible. However, it is permanently
\RR-irreducible which is a contradiction.

It follows that if $(\lambda,\rho)$ is present at the start of the $n$-th
Knuth--Bendix pass, it will be sewn into \WDiff at some point during the $n$-th
Knuth-Bendix pass. As in the proof of Lemma~\ref{stays R-reducible}, the only
way $(\lambda,\rho)$ would not be accepted by $\Rules[n+1]$ is if some proper
prefix or suffix is accepted by $\Rules[n+1]$. But this would contradict
$(\lambda,\rho)$ being a permanent rule. Therefore, $(\lambda,\rho)$ is
accepted by $\Rules[m]$ for each $m\geq n$. 
\end{proof}

\begin{lemma}\label{eventually permanent}
Let $u$ be a fixed string. Then there is a $t_0$ depending on $u$, such that,
for all $t \ge t_0$, each elementary \Storen t-reduction of $u$ is associated
to a permanent rule.
If all proper substrings of $u$ are permanently irreducible, then, for $t\ge t_0$,
there is at most one elementary reduction of $u$, and this is associated
to a permanent rule $(u,w)$. 
\end{lemma}
\begin{proof}
There are only finitely many substrings of $u$. So we need only show
that, given any string $v$, there is a $t_0$ such that for all $t\ge t_0$,
each rule in \Storen t with left-hand side $v$ is permanent.
If there is a proper substring of $v$ which is not permanently irreducible,
then at some time $s_0$ it becomes \Storen{s_0}-reducible.
By Lemma~\ref{stays S-reducible}, it is \Storen s-reducible for $s\ge s_0$.
By Lemma~\ref{stays R-reducible}, it becomes \RR-reducible at the beginning of
the next Knuth--Bendix pass after $s_0$. During this pass all rules
with left-hand side $v$ will be deleted. Also, since this proper substring
of $v$ is now permanently \RR-reducible, no rule with left-hand side
equal to $v$ will ever be inserted subsequently. In this case the lemma
is true since ultimately there are no rules with left-hand side $v$.

So we assume that each proper substring of $v$ is permanently irreducible,
and that $v$ itself is \Store-reducible at some time $t$. A rule $(v,w)$ will 
be permanent if
$w$ is permanently irreducible. Otherwise it will disappear as a result
of minimization and, by Lemma~\ref{no repeats}, never reappear.
There cannot be two permanent rules $(v,w_1)$ and $(v,w_2)$ with $w_1 > w_2$.
For critical pair analysis would produce a new rule $(w_1,w_2)$ during
the next Knuth--Bendix pass, and so $w_1$ would not be permanently
irreducible.
\end{proof}

\begin{figure}
\begin{picture}(200,220)(25,-130)
\put(120,70){$h_0$}
\put(50,50){\framebox(150,5)}
\put(90,50){\framebox(50,5)}
\put(110,50){\framebox(50,5)}
\multiput(90,51)(0,2){2}{\line(1,0){50}}
\multiput(110,52)(0,2){2}{\line(1,0){50}}
\put(110,60){$\lambda_1$}
\put(135,60){$\lambda_2$}
\put(95,40){$u_1^\prime$}
\put(120,40){$u_2^\prime$}
\put(145,40){$u_3^\prime$}

\put(-40,0){\framebox(150,5){}}
\put(0,0){\framebox(50,5){}}
\put(50,0){\framebox(20,5){}}
\multiput(0,1)(0,2){2}{\line(1,0){50}}
\multiput(50,2)(0,2){2}{\line(1,0){20}}
\put(20,10){$\rho_1$}
\put(55,-10){$u_3^\prime$}
\put(110,35){\vector(-4,-1){70}}

\put(140,0){\framebox(150,5){}}
\put(180,0){\framebox(20,5){}}
\put(200,0){\framebox(50,5){}}
\multiput(180,1)(0,2){2}{\line(1,0){20}}
\multiput(200,2)(0,2){2}{\line(1,0){50}}
\put(185,-10){$u_1^\prime$}
\put(220,10){$\rho_2$}
\put(140,35){\vector(4,-1){70}}

%\put(40,30){\dashbox{5}(170,50){}}

\put(-40,-70){\framebox(150,5){}}
\put(0,-70){\framebox(70,5){}}
\multiput(0,-69)(0,1){4}{\line(1,0){70}}
\put(30,-60){$z_1$}
\put(140,-70){\framebox(150,5){}}
\put(180,-70){\framebox(70,5){}}
\multiput(180,-69)(0,1){4}{\line(1,0){70}}
\put(210,-60){$z_2$}

%\put(30,-15){$\wedge$}
\put(30,-45){$\vee$}
\put(33,-45){\line(0,1){35}}
\put(37,-45){$\star$}
\put(7,-45){$\Storen {t^{\prime\prime}}$}

%\put(210,-15){$\wedge$}
\put(210,-45){$\vee$}
\put(213,-45){\line(0,1){35}}
\put(217,-45){$\star$}
\put(187,-45){$\Storen {t^{\prime\prime}}$}

\qbezier(30,-80)(125,-160)(210,-80)
\put(35,-80){\line(-1,0){5}}
\put(30,-85){\line(0,1){5}}
\put(205,-80){\line(1,0){5}}
\put(210,-85){\line(0,1){5}}
\put(203,-78){$\star$}
\put(212,-90){$\Storen {t^{\prime\prime}}$}
\end{picture}
\caption{When the leftward and rightward reductions from $h_0$ are obtained
from rules $(\lambda_1,\rho_1)$ and $(\lambda_2,\rho_2)$ having overlapping
left-hand sides, this diagram shows the time-$t^{\prime\prime}$ Thue
paths that exist after the critical pair associated with the triple
$(u_1^\prime u_2^\prime u_3^\prime,\rho_1 u_3^\prime, u_1^\prime \rho_2)$ has
been resolved to the pair $(z_1,z_2)$.}
\label{non-disjoint left-hand sides}
\end{figure}
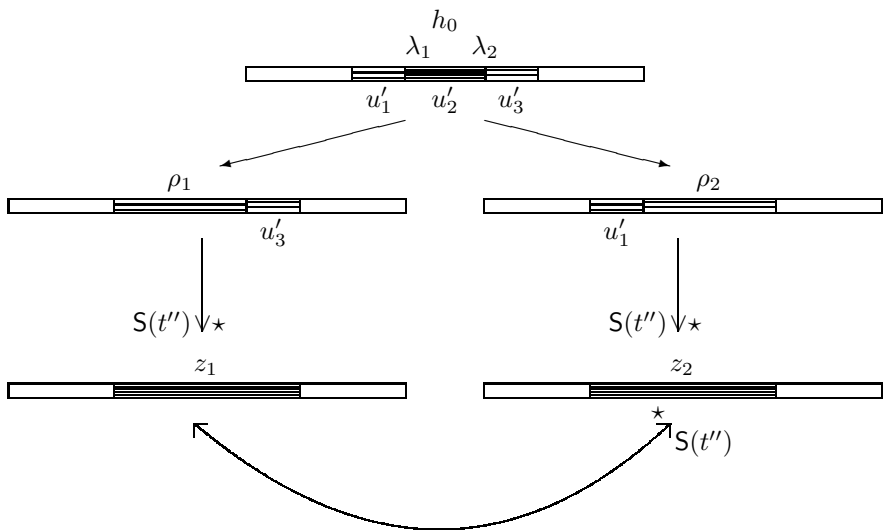

\begin{theorem}\label{confluence}
Let $u$ be a fixed string in $\Astar$ and let $v$ be the smallest
element in its Thue congruence class. Then, for large enough times,
there is a chain of elementary reductions from $u$ to $v$ each associated
to a permanent rule.
After enough time has
elapsed, \RR-reduction of $u$ always gives $v$.
\end{theorem}
\begin{proof}
We start by proving the first assertion.
By hypothesis, we have, for each time $t$, a time-$t$ Thue path
$p_t$ from $u$ to $v$, and we can suppose that $p_t$ contains no repeated
nodes by cutting out part of the path if necessary.
The only reason why we couldn't take $p_{t+1}$ to be $p_t$ is if some
rule $(\lambda,\rho)$, used along the Thue path $p_t$, is deleted at
time $t$. By Lemma~\ref{maintain congruence} we can, however, assume that
each node of $p_{t+1}$ is either already a node of $p_t$ or is smaller
than some node of $p_t$.

Let $h_0$ be the largest node on $p_0$, and suppose that we have already
proved the theorem for all pairs $u$ and $v$ which are connected by a Thue
path with largest node smaller than $h_0$. By induction, using
Proposition~\ref{maintain congruence}, we can assume
that $h_0$ is the largest node on $p_t$ for all time $t$. If $v=h_0$
then since $v$ is the smallest element in its congruence class, there
are no elementary reductions starting from $v$, and we must have $u=v$
in this case.

By Lemma~\ref{eventually permanent}, we may assume that $t_0$ has been
chosen with the property that, for all strings $w\le h_0$ and for all
$t\ge t_0$, all elementary \Storen t-reductions of $w$ are associated
to permanent rules which are accepted by $\Rules[n]$ provided $n$ is
sufficiently large.

Let $h_0 = \mu_t\alpha_t \nu_t\to_{\Storen t} \mu_t\beta_t \nu_t$ be the
rightward elementary reduction of $h_0$ at time $t$.
The rule $(\alpha_t,\beta_t)$ is independent of $t$ for large values
of $t$.
Then $(\alpha_t,\beta_t)$
is permanent and $\alpha_t$ is \RR-reducible for large enough $t$.
If $u\neq h_0$, the same argument applies to
the unique elementary leftward reduction with source $h_0$ at time $t$.

If $h_0 = u$, let $u \rightarrow_{\Storen t} w$ be the first rightward 
elementary reduction for large
values of $t$. By our induction hypothesis, there is a Thue path of
elementary reductions from $w$ to $v$, each associated to a permanent
rule, and with no node larger than $w$, and so we have the required Thue
path from $u$ to $v$. 

Suppose now that $h_0\neq u$, so that we get two permanent rules,
associated to the leftward and rightward elementary reductions of
$h_0$.
If the two elementary reductions are identical, that is, if the two
permanent rules are equal and if their left-hand sides occur in the
same position in $h_0$, then $p_t$ contains a repeated node which we are
assuming not to be the case. So the two elementary reductions occur in
different positions in $h_0$. Now choose $t$ to be large enough so that the two
rules concerned have already been compared in a critical pair analysis
in Step~\ref{compare for overlaps} during some previous $n$-th Knuth--Bendix
pass.

If these two rules have left-hand sides which are
disjoint substrings of $h_0$, then we can interchange their order so as to
obtain a Thue path from $u$ to $v$ where all nodes are strictly smaller than
$h_0$---see \figref{disjoint left-hand sides}. The first assertion of
the theorem then follows by the induction hypotheses in this
particular case.

\begin{figure}
\begin{picture}(200,150)(25,-50)
\put(50,50){\framebox(150,5){}}
\put(80,50){\framebox(30,5){}}
\put(140,50){\framebox(30,5){}}
\multiput(80,51)(0,1){4}{\line(1,0){30}}
\multiput(140,51)(0,1){4}{\line(1,0){30}}
\put(90,60){$\lambda_1$}
\put(150,60){$\lambda_2$}
\put(120,70){$h_0$}
\put(80,40){\vector(-4,-1){70}}
\put(-50,5){\framebox(150,5){}}
\put(-20,5){\framebox(30,5){}}
\put(40,5){\framebox(30,5){}}
\multiput(-20,6)(0,1){4}{\line(1,0){30}}
\multiput(40,6)(0,1){4}{\line(1,0){30}}
\put(-10,15){$\rho_1$}
\put(50,15){$\lambda_2$}
\put(170,40){\vector(4,-1){70}}
\put(150,5){\framebox(150,5){}}
\put(180,5){\framebox(30,5){}}
\put(240,5){\framebox(30,5){}}
\multiput(180,6)(0,1){4}{\line(1,0){30}}
\multiput(240,6)(0,1){4}{\line(1,0){30}}
\put(190,15){$\lambda_1$}
\put(250,15){$\rho_2$}
\put(50,-40){\framebox(150,5){}}
\put(80,-40){\framebox(30,5){}}
\put(140,-40){\framebox(30,5){}}
\multiput(80,-39)(0,1){4}{\line(1,0){30}}
\multiput(140,-39)(0,1){4}{\line(1,0){30}}
\put(70,-5){\vector(4,-1){70}}
\put(180,-5){\vector(-4,-1){70}}
\put(90,-30){$\rho_1$}
\put(150,-30){$\rho_2$}
\put(40,30){\dashbox{5}(170,50){}}
\end{picture}
\caption{Removing the node $h_0$ when the leftward and rightward reductions
are obtained from rules having disjoint left-hand sides.}
\label{disjoint left-hand sides}
\end{figure}
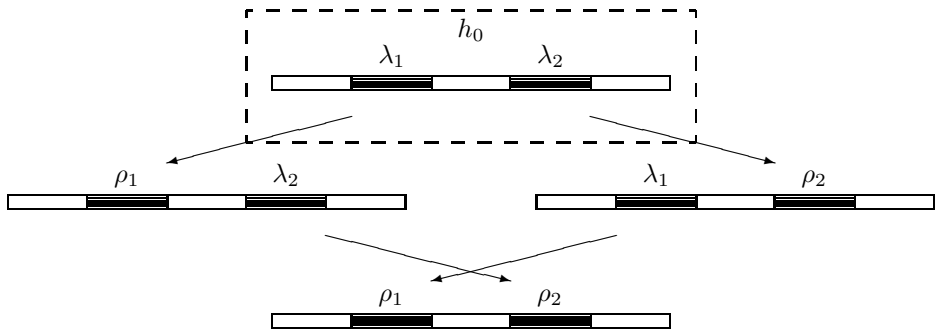

If the two left-hand sides do not correspond to disjoint substrings of
$h_0$ then, by assumption, there is some time
$t^\prime < t$, such that a critical pair
$(u^\prime,v^\prime,w^\prime)$ was considered.
Here
$u^\prime \rightarrow_{\Storen {t^\prime}} v^\prime$ and
$u^\prime \rightarrow_{\Storen {t^\prime}} w^\prime$ are elementary
$\Storen {t^\prime}$-reductions given by the two rules, and $u^\prime$ is a
substring of $h_0$. After the critical pair analysis, at time
$t^{\prime\prime}\leq t$, the Thue paths illustrated in
\figref{non-disjoint left-hand sides} are possible.
As a consequence of
Lemma~\ref{maintain congruence}, it is straighforward to see that for all times
$s\geq t^{\prime\prime}$, $v^\prime$ and $w^\prime$ can be connected by a
time-$s$ Thue path
in which all nodes are no larger than the largest of $v^\prime$ and $w^\prime$.
In particular, this applies at time $t$ so that the targets of the two
elementary
\Storen t-reductions from $h_0$ can be connected by a time-$t$ Thue path in
which all nodes are smaller than $h_0$.
This completes the inductive proof of the first assertion of the
theorem.

We have arranged that $t$ is large enough so that,
for all $w\le u$, all elementary \Storen t-reductions of $w$ are
associated to permanent rules, and such a $w$ can be permanently \RR-reduced
to the least element in its Thue congruence class. It follows that such a $w$ 
is \RR-irreducible if and only if it is minimal in its Thue class. In 
particular \RR-reduction of $u$ must give $v$.
\end{proof}

\begin{corollary} The set of permanent rules in \RR is confluent. The set of
such rules is equal to $\Permanent = \bigcap_{t}\bigcup_{s\ge t}$\Storen s.
A string $u$ is smallest in its Thue congruence class if
and only if it is permanently irreducible.
\end{corollary}
\begin{proof}
The first and third statements are obvious from Theorem~\ref{confluence}.
For the second statement, each permanent rule is contained in \Permanent by
Lemma~\ref{never deleted}. Conversely, if we have a rule $r$ in \Store which is
not permanent, then for all
sufficiently large times $s$ either its right-hand side or a proper substring
of its left-hand side is \Storen s-reducible. Theorem~\ref{confluence} ensures
that this reducible string is \RRn s-reducible for all sufficiently large
times $s$. Therefore $r$ will be minimized and deleted from
\Store. Hence from Lemma~\ref{no repeats} we see that $r$ is not contained in
\Permanent.
\end{proof}

The next result is the main theorem of this paper.
\begin{theorem} Let $G$ be a group with a given finite presentation
and a given ordering of the generators and their inverses.
Suppose that $G$ is \textit{shortlex}-automatic.
Then the procedure given in \ref{main loop} will stabilize at some $n_0$ with
$\Rules[n+1] = \Rules[n]$ if $n\ge n_0$.
\Permanent is then the language of a certain
two-variable finite state automaton and the automaton can be explicitly
constructed.
\end{theorem}
\begin{proof}
\Permanent consists of pairs of strings $(\lambda,\rho)$ giving a valid 
identity in $G$ (by Proposition~\ref{consistent rules}), and where $\rho$
and all proper substrings of $\lambda$ are permanently irreducible.
A string is permanently irreducible if and only if it is the unique
\textit{shortlex} representative of the corresponding group element.
Since this structure is automatic, there are only finitely many word
differences and arrows generated by the rules in
\Permanent. If we therefore weld together the automata $M(r)$ corresponding
to the rules $r$ in \Permanent, we obtain a finite rule automaton $\Rules$.
For each arrow in this automaton, we can pick a specific rule $r$ which
makes use of the arrow when it is read into $\Rules$.

These specific rules will eventually be generated by our Knuth--Bendix
process. Such a rule is never deleted once it is generated, since it is
permanent. So eventually $\Rules[n]$ will contain $\Rules$ as a
sub-automaton. But once this has happened, \RR-irreducible will be
equivalent to \textit{shortlex}-minimal. Therefore all non-permanent
rules will be removed during the next pass, and the redundant states and
arrows of $\Rules[n]$ will be removed. $\Rules[m]$ is then constant for $m>n$.
\end{proof}

Of course, the problem with the above result is that we do not currently
have any method of knowing when we have reached $n_0$. It might be
possible to prove that this question is undecidable if one varies over
all \textit{shortlex} presentations of \textit{shortlex}-automatic groups.
It might also be undecidable in one varies over all finite presentations
of word hyperbolic groups.

\section{Miscellaneous details}
In this section we present a number of points which did not seem to
fit elsewhere in this paper.

\subsection{The structure of \Store}
%We need to be able to insert rules into \Store and delete them.
%Given the left-hand side of a rule $\lambda$, we need to be able
%to find quickly a rule of the form $(\lambda,\rho)$ in \Store, if there is one.
%We also need to be able to choose one of the subsets
%\Considered, \This or \New and iterate through it, obtaining the relevant
%rules one at a time.
This set is given in a data structure arranged so that it is quick to find
a rule in it, given only the left-hand side, quick to delete a specified
rule and quick to add a rule. All these operations take place
repeatedly in the Knuth--Bendix program.
It is also an advantage to have a robust enough method for iterating
through \Store, so that the process is not disrupted if rules are
added or deleted while the iteration is proceeding. (We don't mind if
the iteration fails to catch the rules added during iteration.)

\subsection{Aborting.} It is possible that we come to a situation where
the procedure is not noticing that certain strings are reducible, even though
the necessary information to show that they are reducible is already
in some sense known. It is also possible that reduction is being carried
out inefficiently, with several steps being necessary, whereas in some
sense the necessary information to do the reduction in one step is
already known. An indication that our procedure could be improved is that
\WDiff is constantly changing, with two states being identified and consequent 
welding, or with new states or arrows being added. In this case it might be 
advisable to abort the current Knuth--Bendix pass.

To see if abortion is advisable, we can record statistics about how much
\WDiff has changed since the beginning of a pass. If the changes seem
excessive, then the pass is aborted. A convenient place for the program to
decide to do this is just before another rule from \New is examined
at Step~\ref{process New}.

If an abort is decided upon then all
states and arrows of \WDiff are marked as \needed.
At this point the program jumps to Step~\ref {start main loop}.

\subsection{Priority rules.} A well-known phenomenon found when
using Knuth--Bendix to look for automatic structures, is that rules
associated with finding new word differences or new arrows in
\WDiff should be used more intensively than other rules. Further aspects
of the structure are then found more quickly. These observations are
not the consequence of a theorem---they are observed when programs are
run.

A new rule associated with new word differences or new arrows is
marked as a priority rule. If a priority rule is minimized, the output is
also marked as a priority rule. If a priority rule is added to one of the lists
\Considered, \This or \New, it is added to the front of the list, whereas
rules are normally added to the end of the list.
Just before deciding to add a priority rule to \New, we check to see if
the rule is minimal. If so, we add it to the front of \This instead of to
the front of \New.

When a rule is taken from \This at Step~\ref{process This} during the 
main loop, it is normally compared with all rules in \Considered,
looking for overlaps between left-hand sides. In the case of a priority rule,
we compare left-hand sides not only with rules in \Considered,
but also with all rules in \This. If a normal rule $(\lambda,\rho)$
is taken from \This and comparison with a rule in \Considered gives rise
to a priority rule, then the rule $(\lambda,\rho)$ is also marked as a priority
rule. It is then compared with all rules in \This, once it has been
compared with all rules in \Considered.

Treating some rules as priority rules makes little difference unless there
is a mechanism in place for aborting a Knuth--Bendix pass when \WDiff
has sufficiently changed. If there is such a mechanism, it can make a big
difference.

\subsection{An efficiency consideration.}
During reduction we often know have a state $s$ in a two-variable
automaton. We usually know $x\in A$ and we are looking for an arrow
labelled $(x,y)$ with certain properties, where $y\in A^+$. It therefore
makes a big difference if the arrows with source $s$ are arranged so that
we have rapid access to arrows labelled $(x,y)$ if $x$ is given.

\section{The past and the future}

\subsection{A failed idea.} Our original idea was to avoid having an
explicit finite set of rules \Store. Instead we tried to attach extra information
to the states and arrows of our automata so that the set of rules implicitly
held included both a finite set, corresponding to our current \Store,
and the possibly infinite set held by the automata. The idea was to avoid
using the considerable amount of space used by \Store. This idea did not work
and we now explain why.

The idea was that it didn't matter too much if the finite set of rules held
implicitly was too big. The logic of Knuth--Bendix only goes wrong if it
is too small. However, if the extra information attached to states and arrows
is not sufficiently explicit, there is often a huge
growth from one pass to the next in the finite
set of rules implicitly held. This growth is not caused by the Knuth--Bendix
process itself, but is a by-product of the way we are using the extra
information to specify the finite set of rules.

Another approach is then to attach much more information to states and
arrows in an attempt to limit the unnecessary growth referred to in the 
previous paragraph. But this extra information itself requires a lot of
space, more even than holding the rules separately! Moreover, it turned out
not to be possible to conveniently limit the growth as much as was
necessary. So holding more information in the finite state automata
was worse on all counts, including the complexities of writing the code,
than the simpler scheme of holding the rules separately.

\subsection{The present.}
Many of the ideas in this paper have been implemented
in C++ by the second author. But some of the ideas in this paper only occurred to
us while the paper was being written, and the procedures and algorithms
presented in this paper seem to us to be substantial improvements on what
has been implemented so far. An unfortunate result of this is that we are
unable to present experimental data to back up our ideas, although many
of these ideas have been explored in depth with actual code. Our experimental 
work has been essential in enabling us to come to the
better algorithms which are presented here.

\subsection{Comparison with \textit{kbmag}.} Here we describe the
differences between our ideas and the ideas in Derek Holt's \textit{kbmag}
programs \cite{KBMAG}. These programs try to compute the 
\textit{shortlex}-automatic
structure on a group. Our program is a substitute only for the first program
in the \textit{kbmag} suite of programs.

In \textit{kbmag}, fast reduction is carried out using an automaton with
a state for every prefix of every left-hand side. In our program we also
keep every rule. However, the space required by a single character in our
program is less by a constant multiple than the space required for a
state in a finite state automaton. Moreover, compression techniques
could be used in our situation so that less space is used, whereas
compression is not available in the situation of \textit{kbmag}.

The other large objects in our set-up are the automata $P(\Rules[n])$
defined in Section~\ref{which words are reducible} and $Q(\Rules[n])$
defined in Section~\ref{finding the left-hand side}.
In \textit{kbmag}, there has also to be an automaton like $P(\Rules[n])$, and
it is possible to arrange that this automaton is only constructed
after the Knuth--Bendix process is halted. This can avoid running out of space.
In \textit{kbmag} there is no analogue of our $Q(\Rules[n])$.

In \textit{kbmag}, reduction is carried out extremely rapidly. However,
as new rules are found, the automaton in \textit{kbmag} needs to be
updated, and this is quite time-consuming. In our situation, updating the
automata is quick, but reduction is slower because the string has to be
read into two different automata. Moreover we sometimes need to use the
method of Section~\ref{finding the right-hand side} which is slower
(by a constant factor) than simply reading a string into a deterministic
finite state automaton.

In \textit{kbmag}, there is a heuristic, which seems to be inevitably arbitrary,
for deciding when to
stop the Knuth--Bendix process. In our situation there is a sensible
heuristic, namely we stop if we find $\Rules[n+1] = \Rules[n]$.

In the case of \textit{kbmag}, there are
occasional cases where the process of finding the set of word differences
oscillates indefinitely. This is because redundant rules are sometimes
unavoidably introduced into the set of rules, introducing unnecessary
word differences. Later redundant rules are eliminated and also the
corresponding word differences. This oscillation can continue indefinitely.
Holt has tackled this problem in his programs
by giving the user interactive modes of running them.

In our case,
 the results in Section~\ref{correctness} show that, given a \textit{shortlex}
automatic group, the automaton $\Rules[n]$ will eventually stabilize,
given enough time and space.

We believe that the main advantage of our approach will only become
evident when looking at very large examples. We plan to carry out a
systematic examination of \textit{shortlex}-automatic groups generated
by Jeff Weeks' \textit{SnapPea program}---see \cite{weeks}---in order
to carry out a systematic comparison.

\subsection{Other situations}
We should remark that our methods should apply with some modifications
to certain other orderings, not
only to the \textit{shortlex}-ordering. The essential feature we need
is that the set of pairs $(u,v)$, such that $u>v$, is a regular language.
Other orderings than \textit{shortlex}, for example an ordering called
the \textit{wreath product ordering},
have been useful in theoretical discussions \cite{holthurt}.
The wreath product ordering is used in programs by Holt which look for
coset automatic structures.

Bill Thurston has suggested that we generalize our programs to apply
directly to a triangulated space rather than to a group. It should be
straightforward to make this generalization in both the \textit{kbmag}
programs and in ours.

\bibliography{KB} 
\bibliographystyle{plain}

\end{document}